
\input phyzzx.tex

\def\ex{{\hbox{\rm e}}}

\def\im{{\hbox{\rm Im}}}

\def\tr{{\hbox{\rm Tr}}}
\def\ch{{\hbox{\rm ch}}}

\def\rl{\Lambda_{\hbox{\sevenrm R}}}
\def\wl{\Lambda_{\hbox{\sevenrm W}}}
\def\longl{{\cal L}_{\hbox{\sevenrm R}}}
\def\fc{{\cal F}_{2k+\cox}}
\def\vev{vacuum expectation value}

\def\cox{{g^\vee}}

\tolerance=500000
\overfullrule=0pt
\def\np{Nucl. Phys.}
\def\pl{Phys. Lett.}

\def\pr{Phys. Rev.}

\def\cmp{Comm. Math. Phys.}
\def\ijmp{Int. J. Mod. Phys.}

\def\lmp{Lett. Math. Phys.}
\def\bams{Bull. AMS}
\def\am{Ann. of Math.}
\def\jpsc{J. Phys. Soc. Jap.}
\def\topo{Topology}

\def\knot{Journal of Knot Theory and Its Ramifications}
\def\ams{Ann. Math. Stud.}
\def\mm{Math. Magazine}

\tolerance=500000
\overfullrule=0pt

\tolerance=500000
\overfullrule=0pt

\Pubnum={US-FT/9-95 \cr q-alg/9507031}
\date={July, 1995}
\pubtype={}
\titlepage

\title{A RELATION BETWEEN THE KAUFFMAN
AND THE HOMFLY POLYNOMIALS FOR TORUS KNOTS}
 \author{J.M.F. Labastida\foot{E-mail: LABASTIDA@GAES.USC.ES} and  E. P\'erez}
\address{Departamento de F\'\i sica de Part\'\i culas\break Universidade de
Santiago\break E-15706 Santiago de Compostela, Spain}

\abstract{Polynomial invariants
corresponding to the fundamental representation of the
gauge group $SO(N)$ are computed for arbitrary torus knots
in the framework of Chern-Simons gauge theory making use of knot operators. As
a
result, a formula which relates the Kauffman and  the HOMFLY polynomials
for  torus knots is presented.}

\endpage
\pagenumber=1

\chapter{Introduction}

Knot operators \REF\nos{J.M.F. Labastida
and A.V. Ramallo \journal\pl&B227(89)92,
{\sl Nucl. Phys.} {\bf B} (Proc. Suppl.) {\bf 16B} (1990) 594.}
\REF\llr{J.M.F. Labastida, P.M. Llatas and A.V. Ramallo
\journal\np&B348(91)651.} [\nos,\llr]
have shown to be a powerful tool in Chern-Simons gauge  theory
\REF\witCS{E. Witten \journal\cmp&121(89)351.} [\witCS]
to obtain general expressions for knot
invariants related to torus knots and links.
Computations by other methods
\REF\martin{S.P. Martin \journal\np&B338(90)244.}
\REF\ygw{K. Yamagishi, M.-L. Ge and Y.-S. Wu \journal\lmp&19(90)15.}
\REF\muk{S. Mukhi, ``Skein Relations and Braiding in Topological Gauge
Theory", Tata preprint, TIFR/TH/98-39, June 1989.}
\REF\horne{J.H. Horne \journal\np&B334(90)669.}
\REF\kcp{T.W. Kim, B.H. Cho and S.U. Park\journal\pr&D42(90)4135}
\REF\kg{R.K. Kaul and T.R.
Govindarajan\journal\np&B380(92)293\journal\np&B393(93)392
\journal\np&B402(93)548}
\REF\haya{M. Hayashi\journal\np&B405(93)228}
[\martin,\ygw,\muk,\horne,\kcp,\kg,\haya]
have been succesful for specific knots but
not to obtain general expressions for knot sequences as torus knots.
\REF\poli{J.M. Isidro, J.M.F. Labastida and A.V. Ramallo
\journal\np&B398(93)187.}
Knot operators have
been used in [\poli] where  a formula
for the invariants for torus knots and links
carrying  arbitrary representations of the gauge group $SU(2)$ has been
presented. For the fundamental representation it covers the case of the Jones
polynomial \REF\jones{V.F.R. Jones \journal\bams&12(85)103.}
\REF\jonesAM{V.F.R. Jones\journal\am &  126  (87) 335.}
[\jones,\jonesAM], while for higher dimensional representations it covers
the case of  the Akutsu-Wadati polynomials
\REF\aw{Y. Akutsu and M. Wadati
\journal\jpsc&56(87)839;\journal\jpsc&56(87)3039
Y. Akutsu, T. Deguchi and M. Wadati
\journal\jpsc&56(87)3464;\journal\jpsc&57(88)757; for a review see
M. Wadati, T. Deguchi and Y. Akutsu, Phys. Rep. {\bf 180} (1989)247. }
[\aw]. They have also been used in \REF\lama{J.M.F.
Labastida and M. Mari\~no\journal\ijmp&A10(95)1045.} [\lama]
where a formula for the HOMFLY polynomial
\REF\homfly{P. Freyd, D. Yetter, J. Hoste, W.B.R. Lickorish, K. Millet and
A. Ocneanu \journal\bams&12(85)239.}
[\homfly,\jonesAM] for
arbitrary torus knots and links has been presented.
For the case of torus knots the formula obtained in [\lama]
for the HOMFLY polynomial coincides with the one  presented by Jones in
[\jonesAM], and later reobtained using quantum groups by  Rosso and Jones in
\REF\rosso{M. Rosso and V. Jones\journal\knot&2(93)97}
[\rosso].

Knot operators were constructed in [\nos,\llr] for the gauge group $SU(N)$.
In this paper we will present the form of these operators
for arbitrary simple compact groups.
Then, we will use them to compute
knot invariants for arbitrary torus knots carrying
the fundamental representation of $SO(N)$. As a consequence a formula for the
Kauffman polynomial
\REF\kau{L.H. Kauffman,
{\sl \ams} {\bf 115}, Princeton University Press, 1987}
 [\kau] for this type of knots is obtained. This formula turns out
to be equivalent to the one obtained in
\REF\yoko{Y. Yokota\journal\topo &32 (93)309}
[\yoko] using a different method.
Comparing this formula for the Kauffman polynomial to
the one obtained in [\jonesAM,\rosso,\lama]
for the HOMFLY polynomial we obtain a rather
simple relation between them. Denoting the HOMFLY polynomial for a torus knot
$\{n,m\}$ ($n$ and $m$ are coprime integers, $(n,m)=1$) in terms of its
standard
variables $a$ and $z$ by
$P_{n,m}(a,z)$, and the Kauffman polynomial (Dubrovnik version) by
$Y_{n,m}(a,z)$, we find:
$$
P_{n,m}(a,z)={1\over 2}(Y_{n,m}(a,z)+Y_{n,m}(a,-z)) +
{z\over 2(a-a^{-1})}(Y_{n,m}(a,z)-Y_{n,m}(a,-z)).
\eqn\gaussdos
$$
This is the main new result presented in this paper.
The existence of a formula like
\gaussdos\ is rather remarkable.
In general the Kauffman polynomial contains very many
more terms than the HOMFLY polynomial.
This means that an important cancellation of
terms ocurr in \gaussdos.
Notice also that this formula indicates that at least
for torus knots the Kauffman polynomial
distinguishes more knots than the HOMFLY polynomial.
Two torus knots which have the same
Kauffman polynomial also have the same HOMFLY
polynomial but it might
happen that two torus knots have the same HOMFLY polynomial
but different Kauffman polynomials.
At least for torus knots one can state that the
Kauffman polynomial is more fundamental than the HOMFLY polynomial.

As a byproduct of formula \gaussdos\ it will be obtained in
sect. 4 a formula for the Alexander-Conway
polynomial in terms of the first derivative at $a=1$ of
the corresponding Kauffman polynomial.

The paper is organized as follows.
In sect. 2 we present the generalization of the
construction of knot operators based on
Chern-Simons gauge theory for an arbitrary simple compact gauge group.
In sect. 3 we calculate the
Kauffman polynomial for torus knots obtaining a result
in full agreement with a previous calculation.
In sect. 4. we prove  formula \gaussdos\ and derive a formula
for the Alexander-Conway polynomial for torus knots in terms
of the corresponding Kauffman polynomial.
In sect. 5 we add final comments and remarks on our results.
The conventions used in this paper are conveniently compiled in an appendix.

\endpage
\chapter{Knot Operators for arbitrary simple gauge group}

In this section we present
the generalization of the operator formalism developed in
[\nos,\llr]  for an arbitrary simple compact
gauge group and the construction of the
corresponding knot operators. We begin introducing Chern-Simons gauge theory.
Let $M$ be a boundaryless three-dimensional manifold and let $A$ be a
connection
associated to a principal $G$-bundle for some simple Lie group $G$.
The action  which  defines Chern-Simons gauge theory has the form:
$$
S_k(A) = {k\over 4\pi} \int_M \tr\big(
A \wedge dA + {2 \over 3} A \wedge A \wedge A \big),
\eqn\chern
$$
where $\tr$ is
the trace over the fundamental representation of the simple gauge group
$G$, and, for the moment,
$k$ is an arbitray real number. Under a gauge transformation,
$$
A \rightarrow g^{-1}dg+ g^{-1}Ag,
\eqn\gaugetrans
$$
the action \chern\ transforms as:
$$
S_k(A) \rightarrow S_k(A) - {k\over 12\pi} \int_M \tr\big( g^{-1}dg \wedge
g^{-1}dg \wedge g^{-1}dg \big).
\eqn\stransf
$$
The last quantity is closely
related to the winding number of the map $g:M\rightarrow G$,
which is defined as:
$$
\Upsilon(g) = {1\over 48 \pi} \int_M \epsilon^{\mu\nu\rho} f^{abc} \psi^2
C^a_\mu C^b_\nu C^c_\rho,
\eqn\winding
$$
where $C^a_\mu$ is given by,
$$
g^{-1} \partial_\mu g = C^a_\mu T^a,
\eqn\lastes
$$
being $T^a$, $a=1,\dots,\,{\hbox{\rm dim}} (G)$,
the generators of the simple group $G$.
In \winding\ $f^{abc}$ are
the structure constants corresponding to this group, and
$\psi^2$ the squared length  of the longest simple root of $G$. The quantity
$\Upsilon(g) $ in \winding\ is always $2\pi$ times an integer
\REF\witwin{E. Witten\journal\cmp&92(84)455} [\witwin]. To study its relation
to
the second term on the right
hand side of \stransf\ we must take into account that the
generators can be chosen in such a way that,
$$
\tr(T^a T^b) = - y \psi^2 \delta^{ab},
\eqn\dynkin
$$
where $y$ is the Dynkin index
of the fundamental representation of the simple group
$G$. It is clear from \dynkin\
that this index is independent of the scale chosen for
the gauge group generators. From  \winding\ and \dynkin\ follows that \stransf\
can be written as:
$$
S_k(A) \rightarrow S_k(A) - 2yk \Upsilon(g).
\eqn\stransf
$$
The values of $y$ for $SU(N)$ and $SO(N)$ are
$1/2$ and $1$, respectively. For other groups $y$ is a half-integer
or an integer (see the Appendix).
Therefore, if $k$ is an integer the action \chern\ changes
into $2\pi$ times an integer
and the exponential $\exp(iS_k(M))$ is gauge invariant.
Furthermore, for the case of $SO(N)$
is enough to require $k$ to be half-integer.
Defining:
$$
x=2yk,
\eqn\laequis
$$
one has, in general, the following quantization condition:
$$
x=2yk \in {\bf Z}.
\eqn\quant
$$

For values of $k$ satisfying
the quantization condition \quant\ the partition function
of the theory is defined as,
$$
Z_k(M) = \int [{\cal D}A]_M \exp(i S_k(A)),
\eqn\partfun
$$
where the functional integration is over gauge non-equivalent connections.
This partition function
is a topological invariant because the action $S_k(A)$ does
not depend on the metric on $M$. Other topological-invariant quantities are
constructed introducing
operators in the integrand of the functional integral present
in \partfun. These
operators must be gauge-invariant and metric-independent to lead to
topological-invariant quantities. Wilson lines constitute
an important class of these
operators. Let $\gamma$ be a close curve in $M$ and let $R$ be an irreducible
representation of the gauge group.
The Wilson line operator associated to $\gamma$ and $R$ is:
$$
W_R^\gamma(A) = \tr_R \big( {\hbox{\rm P}} \exp \int_\gamma A \big),
\eqn\wline
$$
where P denotes a path-ordered product along $\gamma$. We will be interested in
computing the vacuum expectation values of products of these operators, \ie,
functional integrations of the form:
$$
\int [{\cal D}A]_M \big(\prod_{i=1}^n W_{R_i}^{\gamma^i} \big) \exp(i S_k(A)).
\eqn\vev
$$

In order to generalize the operator formalism
developped in [\nos,\llr] let us assume
that there are some
Wilson lines $L_i$ on the manifold $M$. We will perform a Heegaard
splitting on $M$ in
such a way that no Wilson line is cut. The case in which this does
not happen has been studied in \REF\bloques{J.M.F. Labastida and A.V.
Ramallo\journal\pl&228(89)214} [\bloques].
In this formalism, the vacuum expectation
values are expressed as an inner
product of states in a Hilbert space. These states are
defined as functional integrals over
configurations on each of the $g$-handlebodies
$M_1$ and $M_2$ which result from the Heegaard splitting.
In order to construct these states let us introduce complex local coordinates
on the Riemman surface $\Sigma$
which corresponds to the common boundary of $M_1$ and
$M_2$,
$$
z=\sigma_1 + i \sigma_2,
\,\,\,\,\,\,\,\,\,\,\,\,\,\,\,\,\,\,\,\,
\bar z=\sigma_1 - i \sigma_2,
\eqn\coord
$$
and let us use complex components
for the part of the gauge connection parallel to the
surface $\Sigma$:
$$
A_z = {1\over 2}(A_1-iA_2),
\,\,\,\,\,\,\,\,\,\,\,\,\,\,\,\,\,\,\,\,
A_{\bar z} = {1\over 2}(A_1+iA_2).
\eqn\lasas
$$
Our aim is to define wave
functionals which will be functional integrals over field
configurations in the
$g$-handlebodies resulting after the splitting with the value of
$A_{\bar z}$ fixed at the boundary. The inner product will be implemented as
an integration over the
components $A_z$ and $A_{\bar z}$ on the common boundary.

Following [\nos,\llr] we will use
the formalism of the holomorphic quantization. Wave
functionals associated to the
$g$-handle body $M_1$ enclosing $p$ Wilson lines are
defined as,
$$
\Psi_1[A_{\bar z}] =
 \int [{\cal D}A]_{M_1} \big(\prod_{i=1}^p W_{R_i}^{\gamma^i} \big)
\exp\big(i S_k(A)-{k\over 2\pi}\int_\Sigma \tr(A_z A_{\bar z})\big),
\eqn\waveuno
$$
where $[{\cal D}A]_{M_1}$ represents the
functional integration measure over gauge
orbits such that $A_{\bar z}$
is fixed at the boundary $\Sigma$. A similar expression
defines the wave functional $\Psi_2[A_{z}]$ for the $g$-handle body $M_2$.
The vacuum expectation value \vev\ is given by the following inner product:
$$
\big(\Psi_2 | \Psi_1 \big) =
\int [{\cal D}A_z {\cal D}A_{\bar z}]_{\Sigma}
\exp\big({k\over \pi}\int_\Sigma \tr(A_z A_{\bar z})\big)
\overline{\Psi_2[A_{\bar z}]}\Psi_1[A_{z}].
\eqn\inner
$$
Let us recall a few important facts
related to this formalism. Boundary terms like the
one in \waveuno\ are introduced to make the wave functional well defined, \ie,
depending on $A_{\bar z}$ on the boundary $\Sigma$. Also, such a term is the
one
responsible for having a functional
integral in \waveuno\ which is extremal for gauge
configurations such that the field strength of $A$
vanishes in the interior of $M_1$.

The commutation relations of the
canonically conjuagte fields $A_z$ and $A_{\bar z}$
on $\Sigma$
can be read from the exponent of
the exponential inserted in \inner. They take the
form:
$$
[A_{\bar z}^a(\sigma), A_z^b(\sigma')] = { \pi \over 2y\psi^2k} \delta^{ab}
\delta^{(2)}(\sigma-\sigma').
\eqn\comre
$$
Our next step is to compute
explicitly the wave functionals \waveuno\ in order to
obtain a description of the Hilbert
space of the theory. To carry this out we will use
standard parametrizations
of the gauge fields $A_z$ and $A_{\bar z}$ on the Riemann
surface $\Sigma$. We will address
the situations corresponding to genus zero and one.

\section{Genus-zero handlebody}

Let $M_1$ be a solid ball and $\Sigma=S^2$ its boundary. On $S^2$ the fields
$A_z$ and $A_{\bar z}$ can be parametrized as:
$$
A_{\bar z} = u^{-1} \partial_{\bar z} u,
\,\,\,\,\,\,\,\,\,\,\,\,\,\,\,\,\,\,\,\,
A_{z} = \bar u^{-1} \partial_{z} \bar u,
\eqn\zero
$$
where $u$ is a single-valued map $u: S^2 \rightarrow G^c$, being $G^c$ the
complexification of $G$. Since $A_{\bar z}^\dagger = - A_z$ one has that
$u^\dagger = \bar u^{-1}$.
The gauge transformations \gaugetrans\ take the following form for fields on
the
surface $S^2$:
$$
A_{\bar z} \rightarrow g^{-1}\partial_{\bar z} g + g^{-1}A_{\bar z} g,
\,\,\,\,\,\,\,\,\,\,\,\,\,\,\,\,\,\,\,\,
A_{z} \rightarrow g^{-1}\partial_{z} g + g^{-1}A_{z} g,
\eqn\gaugetransdos
$$
where $g$ is map $g:S^2\rightarrow G$. In the parametrization \zero\ these
gauge
transformations take the simple form $u\rightarrow ug$.

The next step is to express the measure
$[{\cal D}A_z {\cal D}A_{\bar z}]_{S^2}$ in \inner\ in
terms of an infinite product of
de Haar measures of $G^c$.
This involves the computation of a Jacobian which takes the
form
\REF\polya{A.M. Polyakov and P.B. Wiegmann\journal\pl&B131(83)121}
\REF\gaw{K. Gawedzki and A. Kupiainen\journal\np&B320(89)625}
[\polya,\gaw]:
$$
[{\cal D}A_z {\cal D}A_{\bar z}]_{S^2}
= \exp \big( {g^\vee \over 2y} \Gamma(u\bar
u^{-1}) \big)
|\det \partial_z \partial_{\bar z} | {\hbox{\rm d}} u {\hbox{\rm d}} \bar u,
\eqn\measure
$$
where $g^\vee$ is the dual Coxeter number of $G$ and $\Gamma(\alpha)$ is the
Wess-Zumino-Witten action [\witwin]:
$$
\eqalign{
\Gamma(u) = &{1 \over 2 \pi} \int_{S^2}
\tr ( \alpha^{-1}\partial_z
\alpha \, \alpha^{-1}\partial_{ \overline z} \alpha) \cr &+
{i \over 12 \pi} \int_{M_1} \epsilon^{\mu\nu\rho} \tr (
 {\tilde \alpha}^{-1}\partial_{\mu} {\tilde \alpha} \, {\tilde
\alpha}^{-1}\partial_{\nu}
{\tilde \alpha} \, {\tilde \alpha}^{-1}\partial_{\rho}
{\tilde \alpha} \big).\cr}
\eqn\setas
$$
In \setas\ $\alpha$ is a map
$\alpha : S^2 \rightarrow G$, and $\tilde \alpha$ is one of the
extensions of this map
to the interior of the solid ball $M_1$. The measure
\setas\ does not depend
on the choice of extension of the map $\alpha$. For different
choices, the resulting
Wess-Zumino-Witten actions differ by $2iy$ times an integral of
the form \winding\ where
$M=S^2$. Therefore, since $g^\vee$ is always an integer, the
measure \measure\ is well defined.
It is also important to remark that this measure is
gauge invariant.

In order to write wave
functionals in terms of $u$ and $\bar u$ one would like to
factor the measure \measure\ appropriately.
This is however not obvious due to the
Polyakov-Wiegmann condition
[\polya] satisfied by the Wess-Zumino-Witten action:
$$
\Gamma(\alpha\beta) =
\Gamma(\alpha) + \Gamma(\beta) + \langle \alpha, \beta \rangle,
\eqn\pw
$$
where we have introduced,
$$
\langle \alpha, \beta \rangle = {1\over \pi} \tr (\alpha^{-1}\partial_{\bar z}
\alpha \, \partial_z\beta \, \beta^{-1}).
\eqn\newinner
$$
As in [\nos,\llr],
we will solve this problem making the following choice of measure
on the boundaries of $M_1$ and $M_2$:
take the measure \measure\ without those factors that
only depend on the gauge
variables which are not being integrated over in the path
integral representation of the wave functional.
Working in a gauge where the radial
component of $A$ on $S^2$ vanishes this amounts to choose:
$$
\exp\big( {g^\vee \over 2y}
(\Gamma(\bar u^{-1})+\langle u , \bar u^{-1} \rangle)
\big) {\hbox{\rm d}} \bar u \,\,\,\,\,\, {\hbox{\rm for}}
\,\,\,\,\,\, \Psi_1,
\eqn\coruno
$$
and,
$$
\exp\big( {g^\vee \over 2y} (\Gamma(u)+\langle u , \bar u^{-1} \rangle)
\big) {\hbox{\rm d}}  u \,\,\,\,\,\, {\hbox{\rm for}}
\,\,\,\,\,\, \Psi_2.
\eqn\cordos
$$
In doing this an extra factor $\exp(\langle u , \bar u^{-1} \rangle)$ has been
introduced. One must account for it in \inner. This implies that the
exponential
factor in \inner\ has to be redefined to:
$$
\exp\big({1\over \pi}(k+{g^\vee\over 2y})\int_\Sigma\tr(A_{\bar z} A_z)\big),
\eqn\redefined
$$
so that the inner product \inner\ becomes:
$$
\big(\Psi_2 | \Psi_1 \big) = \int
 {\hbox{\rm d}} u {\hbox{\rm d}} \bar u |\det \partial_z \partial_{\bar z} |
\exp\big({1\over \pi}(k+{g^\vee\over 2y})\int_\Sigma\tr(A_{\bar z} A_z)\big)
\overline{\Psi_2[A_{\bar z}]}\Psi_1[A_{z}],
\eqn\innerdos
$$
where $A_z$ and $A_{\bar z}$ are given by \zero.

As shown in [\nos,\llr],
the form of the wave functional is determined using gauge
invariance. Under the gauge transformations \gaugetransdos, the wave functional
\waveuno\ transforms as:
$$
\Psi[A_{\bar z}]
\rightarrow \Psi[g^{-1} A_{\bar z} g + g^{-1}\partial_{\bar z} g]=
\exp\Big(-(k+{g^\vee\over 2y})\big(
\Gamma(g)+{1\over \pi}\int_\Sigma \tr(A_{\bar z}\partial_z g\,g^{-1})\big)\Big)
\Psi[A_{\bar z}],
\eqn\psitrans
$$
where the variation of
the factor \coruno\ introduced in the measure has been taken
into account.
Notice that in doing the gauge transformation \psitrans\ an extension
to the interior of $M_1$ of the map $g$ on the boundary $\Sigma$ has been done.
The result
\psitrans\ is independent of the choice of extension when $k$ satifies the
quantization condition \quant. The solution to \psitrans\ has the form:
$$
\Psi[A_{\bar z}] = \xi \exp\big(-(k+{g^\vee\over 2y})\Gamma(u)\big).
\eqn\solu
$$
It is known [\witCS] that the Hilbert space  for the case of the solid ball is
one-dimensional.
Independently of the form of the Wilson lines contained in the solid
ball the corresponding wave functional must be proportional to \solu. The wave
functional \solu\ satisfies the Gauss law emanating from the Chern-Simons
action
\chern:
$$
F_{\bar z z}^a \Psi[A_{\bar z}] =0,
\eqn\gauss
$$
where $F_{\bar z z}^a $ are the components of the gauge field strength.
To verify \gauss\ one must use the
commutation relations
for the gauge fields $A_z$ and $A_{\bar z}$ resulting from
\innerdos.

\section{Genus-one handlebody}

In this section we describe
the construction of the operator formalism for the case of
genus one: $\Sigma=T^2$.
The strategy is similar to the one in the previous section. The
non-trivial homology structure
of the torus $T^2$ will provide a richer framework. Let
us first introduce some data to caharacterize the torus.

We will denote the
holomorphic abelian differential of a torus $T^2$ with modular
parameter $\tau$ by $\omega(z)$.  Labeling the homology cycle on $T^2$ which is
contractible in the
handlebody by $\alpha$, and  the one which is not by $\beta$, the
holomorphic form $\omega(z)$ satisfies:
$$
\int_\alpha \omega =1,
\,\,\,\,\,\,\,\,\,\,\,\,\,\,\,
\int_\beta \omega =\tau,
\,\,\,\,\,\,\,\,\,\,\,\,\,\,\,
\int_{T^2} \omega \wedge \bar \omega = \im\tau.
\eqn\omegas
$$

The gauge fields $A_z$
and $A_{\bar z}$ on $T^2$ can be parametrized in the following
way [\gaw]:
$$
A_{\bar z}= (u_a u)^{-1} \partial_{\bar z} (u_a u),
\,\,\,\,\,\,\,\,\,\,\,\,\,\,\,
A_{z}= (u_a \bar u)^{-1} \partial_{\bar z} (u_a \bar u),
\eqn\lasamas
$$
where $u$ is a single-valued map,
$u: T^2 \rightarrow G^c$, and $u_a$ a non-single
valued map,  $u_a: T^2 \rightarrow G$, which takes the form:
$$
u_a = \exp\Big( {i\pi\over \im \tau}\int^{\bar z}\overline \omega(z')
\, a\cdot H - {i\pi\over \im \tau}\int^z \omega(z')\, \bar a \cdot H \Big),
\eqn\laua
$$
where,
$$
a=\sum_{i=1}^l a_i \lambda^{(i)},
\,\,\,\,\,\,\,\,\,\,\,\,\,\,\,
H=\sum_{i=1}^l H_i \lambda^{(i)},
\eqn\lahache
$$
being $\lambda^{(i)}$,
$i=1,\dots,l$, the fundamental weights of a simple group $G$
of rank $l$.
A summary of the group-theoretical
conventions used in this paper is contained in the
Appendix. Notice that $u_a$ is in the maximal torus of $G$ and that
$u^\dagger_a=u_a^{-1}$.
As before, $u^\dagger = \bar u^{-1}$, so that $A_{\bar z}^\dagger = - A_z$.

The generalization of the
measure \measure\ for the case of the torus has the form
[\gaw]:
$$
\eqalign{
[{\cal D}A_z {\cal D}A_{\bar z}]_{T^2} =
& \exp \big( {g^\vee \over 2y} \Gamma(u\bar
u^{-1},C) \big)
|\Pi(a,\tau)|^4 (\im\tau)^l \exp\big(-{g^\vee\over y}\langle u_a,
u_a^{-1}\rangle \big)
\cr &
|\det \partial_z \partial_{\bar z} | {\hbox{\rm d}} u
{\hbox{\rm d}} \bar u  {\hbox{\rm d}} u_a
{\hbox{\rm d}}  u_a^\dagger, \cr}
\eqn\measuremas
$$
where $\Gamma(g,B)$ is the gauged Wess-Zumino-Witten action
\REF\ddp{P. Di Vecchia, B. Durhuus and J.L. Petersen\journal\pl&B144(84)245}
[\ddp],
$$
\Gamma (g,B)=\Gamma (g)-{1\over \pi}\int_\Sigma
\tr\big( g^{-1}B_{\bar z}gB_z  - B_{\bar z}\partial_z g{g}^{-1}
+ {g}^{-1}\partial_{\bar z} gB_z - B_z B_{\bar z} \big),
\eqn\veinte
$$
and
$$
\Pi(a,\tau)=\exp \biggl( {g^{\vee} \psi^2 \pi \over 4{\rm Im} \tau}a^2 \biggr)
\Theta^A_{  g^{\vee},\rho}(a,\tau),
\eqn\lagranpi
$$
being $\Theta^A_{ g^{\vee},\rho}(a,\tau)$
the Weyl antisymmetrized theta function
of level $g^\vee$ (see the Appendix), and,
$$
\rho=\sum_{i=1}^l \lambda^{(i)}.
\eqn\rro
$$
The field $C$ in the measure \measuremas\ is:
$$
C_{\bar z} = u_a^{-1} \partial_{\bar z} u_a,
\,\,\,\,\,\,\,\,\,\,\,\,\,\,\,\,\,\,\,\,
C_{z} = u_a^{-1} \partial_{z} u_a,
\eqn\lasces
$$
while the measure ${\hbox{\rm d}} u_a {\hbox{\rm d}}  u_a^\dagger$
takes the form:
$$
{\hbox{\rm d}} u_a
{\hbox{\rm d}} u_a^\dagger = { d^l a \, d^l \bar a \over (\im\tau)^l}.
\eqn\haar
$$

The measure \measuremas\ is invariant under the gauge transformations
\gaugetransdos\ which now take the form,
$$
u\rightarrow ug,
\eqn\filo
$$
which will be called of type (i); under transformations which leave the
fields $A_z$ and $A_{\bar z}$ invariant,
$$
u\rightarrow \hat g^{-1}u,
\,\,\,\,\,\,\,\,\,\,\,\,\,\,\,
u_a \rightarrow u_a\hat g,
\eqn\mena
$$
where $\hat g$ is a map from $T^2$ into the Cartan torus of $G$;
and under modular transformations. This last set of transformations
is described in the Appendix.
The transformations \mena, which will be denoted as type
(ii), involve maps $\hat g$ which  are labelled in the following way:
$$
\hat g_{m,n} = \exp\Big({2\pi i \over \psi^2\im\tau}\big(
(n+m\tau)\cdot H \int^{\bar z}\overline{\omega(z')} - (n+m\bar\tau)\cdot
H\int^z\omega(z)\big)\Big),
\eqn\raices
$$
being $n$ and $m$
elements of the lattice generated by the long roots of $G$, which
will be denoted
by $\longl$, \ie, $n,m\in\longl$. Notice that the maps \raices\ are not
connected to the identity map.

The analogue of the Polyakov-Wiegmann condition \pw\ for the case of the
Wess-Zumino-Witten action [\witwin] takes the form:
$$
\Gamma(u {\overline u}^{-1},C)= \Gamma (u) + \langle u_a,u \rangle +
 \Gamma( \overline u^{-1})
 + \langle \overline u^{-1},\overline u_a^{-1} \rangle -
\langle u_a,u_a^{-1} \rangle
 - {1 \over \pi} \int_{\Sigma}\tr(C_z C_{\overline z}).
\eqn\propie
$$
This expression leads to similar  factorization problems as the ones found
from
\pw. Following [\nos,\llr]  we take
$$
\exp\big({g^\vee\over 2y}(\Gamma(u\bar u^{-1},C)-\Gamma(u)-\langle
u_a,u\rangle)
\big) {\hbox{\rm d}}\bar u {\hbox{\rm d}} u_a^\dagger
\,\,\,\,\,\,\,   {\hbox{\rm for}}   \,\,\,\,\,\,\,  \Psi(A_{\bar z}),
\eqn\une
$$
and,
$$
\exp\big({g^\vee\over 2y}(\Gamma(u\bar u^{-1},C)-\Gamma(\bar u^{-1})-\langle
\bar u^{-1},u_a^{-1}\rangle) \big)
{\hbox{\rm d}} u {\hbox{\rm d}} u_a
\,\,\,\,\,\,\,   {\hbox{\rm for}}   \,\,\,\,\,\,\,  \Psi(A_{z}).
\eqn\dos
$$
After comparing the products of these
two factors to the one in \measuremas\ one finds
that the inner product \inner\ now takes the form:
$$
\eqalign{
\big(\Psi_2 | \Psi_1 \big)
= & \int {\hbox{\rm d}}u \, {\hbox{\rm d}} \overline u \,
{\hbox{\rm d}}u_a \, {\hbox{\rm d}} \overline u_a
|\Pi(a,\tau)|^4  (\im\tau)^{l} \exp \big(-{g^\vee\over 2y}
\langle u_a , u_a^{-1} \rangle \big)   \cr
& \times \exp \bigg({1 \over {\pi }}(k+{g^{\vee}\over 2y}) \int_{\Sigma}
\tr(A_z A_{\overline z}) \bigg)
\overline {\psi_2[A_{\bar z}]} \psi_1[A_{\bar z}].
\cr}
\eqn\gafas
$$

As in the genus-zero case,
the general form of the wave-functional is obtained using
arguments based on its  properties under symmetry transformations.
Performing a gauge transformation of type (i) \filo\ one finds:
$$
\Psi[A_{\overline z}]
\rightarrow \exp \big(-(k + {g^{\vee}\over 2y}) (\Gamma (g)
+\langle u_a u,g \rangle ) \big) \Psi[A_{\overline z}].
\eqn\tipouno
$$
Using the Polyakov--Wiegmann
condition \pw\ one finds that the solution to \tipouno\
can be written as:
$$
\Psi[A_{\overline z}] = \xi \psi_{2yk+{g^\vee}}(u_a u)
\Lambda(u_a),
\eqn\paz
$$
where $\xi$ is a constant, $\Lambda(u_a)$ is arbitrary, and
$\psi_{2yk+{g^\vee}}(u_a u)$ is a functional which satisfies:
$$
\psi_{2yr}(u_a v) = \psi_{2yr}(u_a)
\exp\big(-r(\Gamma(v) + \langle u_a,v\rangle)\big),
\eqn\andrade
$$
for any single-valued map $v : T^2\rightarrow G$.

To search for
solutions to \andrade\ let us perform a symmetry transformation of type (ii)
\mena. One finds,
$$
\Psi[A_{\overline z}] \rightarrow \exp \big({g^{\vee}\over 2y} (\Gamma (\hat g)
+\langle u_a,\hat g \rangle ) \big) \Psi[A_{\overline z}],
\eqn\ada
$$
which implies the following property for $\Lambda(u_a)$ in  \paz:
$$
\Lambda(u_a\hat g) = \exp\big( {g^\vee\over 2y}(\Gamma(\hat g)+
\langle u_a,\hat g \rangle )\big) \Lambda(u_a).
\eqn\romulo
$$
Comparing to \andrade\ it turns out that $\Lambda(u_a)$ and $\psi_{2yr}(u_a)$
are related in the following way:
$$
\Lambda(u_a) = \big[\psi_{{g^\vee}}(u_a) \big]^{-1}.
\eqn\remo
$$

We need now to solve for \andrade. Let us consider the situation in which
the map $v$ is a map as in \raices\ of the form $\hat g_{n_{[i]},0}$ with
$n_{[i]}=\sum_{j=1}^l n_{[i]}^j\alpha_{(j)}$ and $n_{[i]}^j=\delta_i^j$,
being $\alpha_{(j)}$ the simple roots of the group $G$.
Equation \andrade\ takes the form:
$$
\psi_{2yr}(u_{a+n_{[i]}}) = \exp\big(2yr({\pi\over\psi^2\im\tau}
n_{[i]}\cdot n_{[i]} + {\pi\over\im\tau} a\cdot n_{[i]})\big) \psi_{2yr}(u_a).
\eqn\cesar
$$
For maps of the form  $g_{0,m_{[i]}}$ with
$m_{[i]}=\sum_{j=1}^l m_{[i]}^j\alpha_{(j)}$ and $m_{[i]}^j=\delta_i^j$,
one finds,
$$
\psi_{2yr}(u_{a+m_{[i]}\tau}) =
\exp\big({2yr}({\pi\tau\bar\tau\over\psi^2\im\tau}
m_{[i]}\cdot m_{[i]} +
{\pi\over\im\tau} a\cdot m_{[i]}\bar\tau)\big) \psi_{2yr}(u_a).
\eqn\neron
$$
The two types of maps
under consideration generate the maps \raices\ as described in
[\llr]. The general solution
to equations \cesar\ and \neron\ can be expressed in
terms of theta functions of level $r$:
$$
\psi_{2yr,p}(a,\tau) = \exp \big( {y r \pi \psi^2 a^2 \over 2\im\tau}\big)
\Theta_{2yr,p}(a,\tau),
\eqn\herodes
$$
where $p$ is an element
of the weight lattice modulo $2yr$ times the root lattice,
\ie, $p\in \wl / 2yr \rl$. The properties of the theta functions
$\Theta_{2yr,p}(a,\tau)$ are briefly
summarized in the Appendix.

Our analysis leads to the following form for the wave functional:
$$
\Psi[A_{\bar z}] = \xi \exp\big(-(k+{g^\vee\over 2y})(\Gamma(u)+
\langle u_a , u \rangle ) \big){\psi_{2yk+{g^\vee}}(u_a)\over
\psi_{{g^\vee}}(u_a)},
\eqn\cayo
$$
where $\xi$ is a constant, and $\psi_{2yk+{g^\vee}}(u_a)$ and
$\psi_{{g^\vee}}(u_a)$ represent certain linear combinations of
the solutions \herodes.
As shown in [\nos,\llr] the $u$-dependence of the wave functional can be
integrated out obtaining an effective theory. Using \cayo, the inner product
\gafas\ becomes:
$$
\eqalign{
\big(\Psi' | \Psi \big) = & \int
{\hbox{\rm d}}u_a \, {\hbox{\rm d}} u_a^\dagger
|\Pi(a,\tau)|^4  (\im\tau)^{l} \exp \big(-(k+{g^\vee\over y})
\langle u_a , u_a^{-1} \rangle \big)   \cr
& \times
\bar\xi'\xi\overline{ \Bigg[{\psi'_{2yk+{g^\vee}}(u_a)
\over \psi'_{{g^\vee}}(u_a)}
\Bigg] } {\psi_{2yk+{g^\vee}}(u_a) \over \psi_{{g^\vee}}(u_a)}
\int {\hbox{\rm d}}u \, {\hbox{\rm d}} \overline u
\exp\big(-(k+{g^\vee\over 2y})\Gamma(u\bar u^{-1},C)\big),\cr}
\eqn\graco
$$
which, after using the result [\gaw]:
$$
\int {\hbox{\rm d}}u \, {\hbox{\rm d}} \overline u
\exp\big(-(k+{g^\vee\over 2y})\Gamma(u\bar u^{-1},C)\big) =
(\im\tau)^{-{l\over 2}} | \Pi(a,\tau) |^{-2}
\exp({g^\vee\over 2y} \langle u_a, u_a^{-1} \rangle),
\eqn\tiberio
$$
becomes:
$$
\eqalign{
\big(\Psi' | \Psi \big) = & \int
{\hbox{\rm d}}u_a \, {\hbox{\rm d}} u_a^\dagger
|\Pi(a,\tau)|^2  (\im\tau)^{l\over 2} \exp \big(-(k+{g^\vee\over y})
\langle u_a , u_a^{-1} \rangle \big)   \cr
& \,\,\,\,\,\,\,\,\,\,\,\,\,\,  \times
\bar\xi'\xi\overline{ \Bigg[{\psi'_{2yk+{g^\vee}}
(u_a) \over \psi'_{{g^\vee}}(u_a)} \Bigg] }
{\psi_{2yk+{g^\vee}}(u_a) \over \psi_{{g^\vee}}(u_a)}. \cr}
\eqn\graconcio
$$

Weyl invariance forces to choose antisymmetric combinations of the solutions
\herodes. Defining:
$$
\lambda_{2yr,p}(a,\tau) = \sum_{w\in W} \epsilon(w)\psi_{2yr,w(p)}(a,\tau),
\eqn\marcoantonio
$$
where $W$ is the
Weyl group and $\epsilon(w)$  the signature of the element $w\in W$,
the effective inner product \graconcio\ becomes:
$$
\eqalign{
\big(\lambda_{2yk+{g^\vee},q}| \lambda_{2yk+{g^\vee},p}\big) & =
|\xi|^2 \int
{\hbox{\rm d}}^l a \, {\hbox{\rm d}}^l \bar a (\im\tau)^{-{l\over 2}}
\, \exp\big(-(2yk+{g^\vee}){\pi\psi^2\over 2\im\tau} a\cdot\bar a\big)
\cr &
\,\,\,\,\,\,\,\,\,\,\,\,\,\,\,\,\,\,\,
\times
\overline{\lambda_{2yk+{g^\vee},q}(a,\tau)}
\lambda_{2yk+{g^\vee},p}(a,\tau).\cr}
\eqn\cleopatra
$$
{}From this inner product for the effective theory one can read the commutation
relations of its basic operators:
$$
[\bar a^i, a_j] = {2\im\tau \over \pi(2yk+g^\vee)\psi^2}\delta^i_j.
\eqn\comrel
$$
The states $\lambda_{2yr,p}$ of the form \marcoantonio\ which are
independent in
$\wl/2yr\rl$ constitute the physical states  or Hilbert space of
the theory. The set of weights labeling those states constitute the fundamental
chamber ${\cal F}_{2yr}$.

Knot operators are associated to Wilson lines. They correspond to
the form of these operators when represented in the framework of
the Hilbert space which has been constructed. Let us consider a torus
knot labelled by two coprime integers $n$ and $m$, and their corresponding
Wilson line:
$$
W_\Lambda^{(n,m)} = \tr_\Lambda\big( {\hbox{\rm P}} \exp \int_{n,m} A\big).
\eqn\willine
$$
We use the convention in which $n$ ($m$) denotes the number of times
that the Wilson line winds along the $\beta$-cycle ($\alpha$-cycle)
on the torus.

We are interested in the form of this operator when the single valued map
$u$ in  \lasamas\ has been integrated out. In other words, we need the
experssion for the Wilson line \willine\ when $u=1$. Using \laua\ it turns
out to be:
$$
\eqalign{
W_\Lambda^{(n,m)} &= \tr_\Lambda\bigg( \exp\Big( {i\pi\over \im\tau}
\big( (n\bar\tau + m) a\cdot H - (n\tau+m)\bar a\cdot H)\big)\Big)\bigg)
\cr & =
\sum_{\mu\in M_\Lambda} \exp\Big( - {\pi \over \im\tau}(n\bar\tau+m)
a\cdot \mu + {2(n\tau + m)\over
(2yk+g^\vee)\psi^2} \mu\cdot{\partial\over \partial a}\Big), \cr}
\eqn\opeuno
$$
where in the last
step we have used \comrel, and the fact that $H$ is made out of
diagonal matrices whose entries are related to
the components of the set of weights $\mu\in
M_{\Lambda}$,
being $M_{\Lambda}$ the set of weights corresponding to an irreducible
representation of
highest weight $\Lambda$. Using the standard properties of the theta
functions which are compiled in the Appendix one finds:
$$
W_\Lambda^{(n,m)} \lambda_{2yk+{g^\vee},p} =
\sum_{\mu\in M_\Lambda} \exp\Big({2 i \pi \mu^2 n m \over \psi^2(2yk+g^\vee)}
+{4 i \pi m \, p\cdot\mu \over \psi^2(2yk+g^\vee) } \Big)
\lambda_{2yk+{g^\vee},p+n\mu}.
\eqn\venator
$$
These operators are called knot operators. They satisfy the following important
relation:
$$
W_\Lambda^{(1,0)}|\rho\rangle = |\rho+\Lambda\rangle,
\eqn\creator
$$
where $|\rho\rangle$ is
the state corresponding to the weight \rro. As discussed in
[\nos,\llr],
this relation allows to think of the operators $W_\Lambda^{(1,0)}$ as
creation operators
since they create the state corresponding to the highest weight
$\Lambda$ when acting on the vacuum state $|\rho\rangle$.

One important  ingredient in the computation of knot invariants for torus knots
is the knowledge of the corresponding
representation on the set of homeomorphisms
on $T^2$. These homeomorphisms are
generated by modular transformations $S$ and $T$ on
$T^2$ which possess the following
representation \REF\kac{V. Kac, {\it Infinite dimensional Lie algebras},
Cambridge University Press, 1983} [\kac]:
$$
\eqalign{
T_{p,p'} &= \delta_{p,p'}\ex^{2\pi i (h_p - {c\over 24})}, \cr
S_{p,p'} &= {i^{|\Delta_+|}\over (2yk + g^\vee)^{l\over 2}}
\Big({ {\hbox{\rm Vol}}\, \longl^* \over {\hbox{\rm Vol}}\, \longl} \Big)
\sum_{w\in W} \epsilon(w)\ex^{-{4\pi i p\cdot w(p') \over \psi^2
(2yk+g^\vee)}},\cr}
\eqn\modu
$$
where $|\Delta_+|$ is the number of positive roots, $\longl$
is the lattice of long
roots and $\longl^*$ its dual. In \modu\ $h_p$ and $c$ represent the comformal
weight and central charge of the
corresponding two-dimensional conformal field theory:
$$
h_p={p^2-\rho^2 \over \psi^2(2yk+g^\vee)},
\,\,\,\,\,\,\,\,\,\,\,\,
c = {2yk \, {\hbox{\rm dim}}(G) \over 2yk + g^\vee}.
\eqn\cyh
$$

Knot operators provide
a very useful tool to compute knot invariants in lens spaces.
These spaces are
boundaryless three-dimensional manifolds which can be built by joint
of two tori.
The gluing is carried out by an homeomorphism whose representation in the
Hilbert space  which we have
constructed is written in terms of the generators \modu.
If we denote this
representation by $F$, the vacuum expectation value for a Wilson
line corresponding to a torus
knot carrying an irreducible  representation  of highest
weight $\Lambda$ of a simple group $G$ is:
$$
V_\Lambda^{(n,m)}\Big|_F = {\langle\rho|F W_\Lambda^{(n,m)}|\rho\rangle
                      \over \langle \rho | F | \rho \rangle}.
\eqn\vacio
$$
To connect with the standard form in which polynomial invariants are written
we need to correct
\vacio\ in three aspects. Fisrt of all in \vacio\ a choice of frame
for the knot and the manifold has been done. Invariants are usually
expressed in the standard
frames and we must correct \vacio\ so that the contribution
from the knot framing
factor is cancelled, and that the appropiate choice of $F$ is
made. Taking the three-sphere
as our choice of lens space, which will be the case of
interest in this paper,
the standard frame is accomplished considering $F=S$, being $S$ one of the
two generators of modular transformations.
As shown in [\llr] the correction relative to the frame of the knot  is easily
accomplished multiplying by
$$
\ex^{-2\pi i n m h_{\rho+\Lambda}},
\eqn\framing
$$
where $h_{\rho+\Lambda}$ is the conformal weight given in \cyh.
The second aspect
leading to an additional correction for \vacio\ is the fact that the
orientation
chosen for the torus $T^2$ is the opposite
to the standard one. We must therefore do
the following change $m\rightarrow -m$.
Finally, the third aspect is that usually
knot invariants are
normalized in such a way that their value for the unknot
is one. We must
therefore normalize \vacio\ by its value for the unknot.
These three aspects lead to the following proposition:

{\bf Proposition 2.1}. The normalized
knot invariant for a torus knot $\{n,m\}$ in the
standard framing, carrying a $G$ irreducible representation of highest weight
$\Lambda$ on $S^3$ in the standard framing, is:
$$
\eqalign{
X_\Lambda^{(n,m)} & = \ex^{2\pi i n m h_{\rho+\Lambda}}
{V_\Lambda^{(n,-m)} |_{S^3} \over V_\Lambda^{(1,0)} |_{S^3} } \cr
& = \ex^{2\pi i n m h_{\rho+\Lambda}} {\langle\rho | S
W_\Lambda^{(n,-m)} | \rho\rangle
\over \langle\rho |S W_\Lambda^{(1,0)} | \rho\rangle}. \cr}
\eqn\equis
$$
The structures of the knot operators \venator\ and the matrix $S_{p,p'}$
 in \modu\ allow to express this invariant
in terms of the variable
$$
t =\ex^{2\pi i \over 2yk+g^\vee},
\eqn\late
$$
which encloses all
the dependence on $k$. The main purpose of this paper is to compute
\equis\ for the
fundamental representation of the group $SO(N)$. This will lead to the
Kauffman polynomial
[\kau] for torus knots. The resulting formula agrees with the one
given in [\yoko].
The comparison of this
formula to the corresponding known expression for the HOMFLY polynomial
[\jonesAM,\rosso,\lama] will allow to prove \gaussdos.

\endpage
\chapter{Kauffman polynomial for torus knots}

In this section we will make
use of proposition 2.1 to compute the Kauffman polynomial
for torus knots.
We must evaluate \equis\ for the fundamental representation of
$SO(N)$, \ie, we must make $\Lambda=\lambda^{(1)}$. The result is stated in the
following theorem:

{\bf Theorem 3.1.} The Kauffman polynomial for a torus knot $\{n,m\}$ is given
by:
$$\eqalign{
X_{\lambda^{(1)}}^{(n,m)} = &{{ {[1]}\,\lambda^{nm}} \over {
{{[1] + [ 0;1 ] } } }}  \Biggl(
 \sum\limits_{\gamma +\beta +1=n \atop \beta, \gamma \ge 0} t^{-
{m \over 2}(\beta - \gamma)}
\lambda^{-m}
   (-1)^{\gamma}
  \bigg(  {1 \over
{ [n]}} + { 1 \over {  [\beta - \gamma ; 1] }} \bigg) \cr
&\,\,\,\,\,\,\,\,\,\,\,\,\,\,\,\,\,\,\,\,\,\,
\times {1 \over {[\beta]!\,\ [\gamma]!}}
 \prod\limits_{j=-\gamma }^{\beta }
[ j; 1] + \cases{ 0, &n odd; \cr 1, &n even; \cr}
\;\ \Biggr) \cr}
\eqn\viola
$$
where:
$$
[p]=t^{p \over 2}-t^{-{p \over 2}},
\,\,\,\,\,\,\,\,\,\,\,
[p;y] = t^{ p \over 2}\lambda^{y} -  t^{- {p \over 2}} \lambda^{-y},
\,\,\,\,\,\,\,\,\,\,\,
\lambda= t^{N-1 \over 2},
\,\,\,\,\,\,\,\,\,\,\,
t=\ex^{2\pi i \over 2k+ \cox},
\eqn\limon
$$
with $g^\vee = N-2$.

{\bf Proof.} The rest of this
section deals with the proof of this theorem. As $SO(N)$
has two different algebras,
depending on whether $N$ is odd or even, we will have to
study both cases separately.
We will begin with $SO(2l+1)$, $B_l$ being the corresponding
algebra. The main feature of this case is that   the simple roots
of $B_l$ are not all of the same length.
 Notice that since an $\{ n,m\}$
torus knot is isotopically equivalent to the
$\{ -n,-m \}$
torus knot, we can restrict ourselves to torus knots with $n>0$.
Also we will consider the case in which $l>n$. Our results,
however, as in the case of the HOMFLY polynomial computed in [\lama], are
 valid
for arbitrary $l$.
In this proof we make the following choice of normalization for the long roots:
$$
\psi^2 = 2.
\eqn\normaroot
$$
Notice also that for $SO(N)$ the Dynkin index for the fundamental
representation is $y=1$ and therefore \late\ becomes:
$$
t= \ex^{2\pi i \over 2k+g^\vee}.
\eqn\malena
$$

\section{$SO(2l+1)$}

Let us begin
working out the action of the knot
operator  $W^{(n,m)}_{\lambda^{(1)}}$ on the vacuum
state. Using \venator, \normaroot, and the form  of $t$ in \malena, we have:
$$
 W^{(n,-m)}_{\lambda^{(1)}}|\rho\rangle= \sum_{i=1}^{2l+1}
t^{- {1 \over 2}\mu_i^2
 nm - m \mu_i\cdot \rho} | \rho + n \mu_i \rangle
\eqn\paros
$$
where $\mu_i$, $i=1,...,2l+1$, are the weights in $M_{\lambda^{(1)}}$ whose
explicit expression is given in (A.20).
Following the framework
described in the previous section, we must find the canonical
representatives in the fundamental chamber $\fc$
(notice that $2yr=2yk+g^\vee$ and $y=1$) of the weights appearing in
the sum.
The weights present in \paros\ have the following structure:
$$\eqalign{
&\rho+n\mu_1 = (n+1,1,\dots,1), \cr
 &\,\,\,\,\,\,\,\,\,\,\, \vdots, \cr
&\rho+n\mu_j =  (1,\dots,1,1-n,\buildrel j \over {1+n},1,\dots,1),\cr
&\,\,\,\,\,\,\,\,\,\,\,  \vdots, \cr
 &\rho+n\mu_l = (1,\dots,1,1-n,1+2n), \cr
&\rho+n\mu_{l+1} = \rho, \cr
& \rho+n\mu_{l+2} = (1,\dots,1,1+n,1-2n), \cr
& \,\,\,\,\,\,\,\,\,\,\,  \vdots, \cr
& \rho+n\mu_{l+1+j}
=  (1,\dots,1,1+n,\buildrel l-j+1 \over {1-n},1,\dots,1), \cr
 & \,\,\,\,\,\,\,\,\,\,\, \vdots, \cr
& \rho+n \mu_{2l+1} = (1-n,1,\dots,1). \cr}
\eqn\apple
$$
Every weight in the weight lattice can be written as
$w(\mu)+(2k+\cox)\alpha$,
where $w$ is an element of the Weyl group, $\alpha$ a long
root, and $\mu$ is a weight whose components are non-negative. In
the Hilbert space constructed in the previous section the weigths
which possess one or more components which vanish are represented by null
vectors. Since $2l+1>n$
there is no need to add terms of the form $(2k+\cox)\alpha$ to
the weigths in
\apple\ to bring them to a form in which their components are non-negative.
A series of Weyl reflections will be sufficient.
If $n=1$ all the weigths in \apple\ except the first one and $\rho+n\mu_{l+1}$
have one vanishing
component and therefore there are only these two contributions
in the sum present in
\paros. If $n>1$, notice first that the  weights
 $\rho+n\mu_1$ and $\rho+n\mu_{l+1}$ in \apple\ are already in $\fc$.
 For the rest we have the
following cases:

1. Case $i=2, \dots, l$:

a) $2\le i\le n$. We perform the chain of Weyl reflections:
$$
\eqalign{
&\rho + n \mu_i \buildrel \sigma_1 \over \longrightarrow \dots
\buildrel \sigma_{i-2} \over \longrightarrow
 \buildrel \sigma_{i-1} \over \longrightarrow \nu_i = ( n+ 1- i,1,\dots,
\buildrel i \over 2,\dots,1), \;\; i= 2,\dots, l-1, \cr
& \,\ \cr
&\rho + n \mu_l \buildrel \sigma_1 \over \longrightarrow \dots
\buildrel \sigma_{l-2} \over \longrightarrow
 \buildrel \sigma_{l-1}
\over \longrightarrow \nu_l = ( n+ 1- l,1,\dots,1,3). \cr}
\eqn\casandra
$$

The weight $i=l$ will not
be considered as we restrict ourselves to $n<l$.

b) $i>n$. The chain of Weyl reflections is like the one in \casandra:
$$
\eqalign{
&\rho + n \mu_i \buildrel \sigma_{i+1-n} \over \longrightarrow \dots
 \buildrel \sigma_{i-1} \over \longrightarrow  = (1,\dots,1,
\buildrel {i-n} \over 0,1,\dots,\buildrel i \over 2,\dots,1), \;\;\;
 i= 1,\dots, l-1, \cr
& \,\ \cr
 & \rho + n \mu_l \buildrel \sigma_{l+1-n} \over \longrightarrow \dots
 \buildrel \sigma_{l-1} \over \longrightarrow  = (1,\dots,1,
\buildrel {l-n} \over 0,\dots,1,3). \cr}
\eqn\eliane
$$
After $n+1$ reflections the weights get a vanishing component and therefore all
these weights correspond to null vectors and do not contribute to the sum in
\paros. This fact is very important in this calculation because it implies
that the sum \paros\ is truncated. Its upper limit turns out to be $n$
instead of $2l+1$.

2. Case $i=l+2, \dots, 2l+1$:

As $i>n$ for the weights in this case,
we would expect that all of them would achieve
 a vanishing component after a chain
of Weyl reflections. What actually happens is that
for $n$ odd an extra weight will contribute:
$$
\rho + n \mu_{l+2}
 \buildrel \sigma_l \over \longrightarrow  = (1,\dots,1,2-n,2n-1)
 = \rho \;\;\; {\rm for}\;\; n=1. $$
For $j = 2, \dots, l$, one has the following situations:
$$
\eqalign{
 &n \leq j,\cr
& \rho + n \mu_{l+1+j}
 \buildrel \sigma_{l-j+n-1} \over \longrightarrow \dots
 \buildrel \sigma_{l-j+1} \over \longrightarrow
\buildrel \sigma_{l-j+1} \over \longrightarrow = (1,\dots,1,
\buildrel {l-j} \over 2,\dots,1, \dots,
 \buildrel {l-j+n} \over 0,\dots,1); \cr
& \,\ \cr
 &j< n < 2j-1,\cr
& \rho + n \mu_{l+1+j} \buildrel \sigma_{l+j-n} \over \longrightarrow\dots
 \buildrel \sigma_l \over \longrightarrow \buildrel \sigma_{l-1} \over
\longrightarrow
\dots \buildrel \sigma_{l-j+1} \over \longrightarrow = (1,\dots,1,
\buildrel {l-j} \over 2,\dots,1 \dots
 \buildrel {l+j-1-n} \over 0,\dots,1); \cr
& \,\ \cr
& n = 2j-1, \cr
&\rho + n \mu_{l+1+j}  \buildrel \sigma_{l-j+1} \over \longrightarrow \dots
 \buildrel \sigma_l \over \longrightarrow \buildrel \sigma_{l-1} \over
\longrightarrow
 \dots \buildrel \sigma_{l-j+1} \over \longrightarrow = (1,\dots,1,2j-n,
n-(2j-2) \dots,1) = \rho; \cr
& \,\ \cr
& 2j-1 < n <l,\cr
&\rho + n \mu_{l+1+j}  \buildrel \sigma_{l+j-n+1} \over \longrightarrow \dots
 \buildrel \sigma_l \over \longrightarrow \buildrel \sigma_{l-1} \over
\longrightarrow
\dots \buildrel \sigma_{l-j+1} \over \longrightarrow = (1,\dots,1,
\buildrel {l+j-n} \over 0,\dots,1 \dots
 \buildrel {l-j+1} \over 2,\dots,1). \cr }
\eqn\rita
$$
 We see that all the vectors have a vanishing component except when $n=2j-1$,
where the weight $\nu_{l+1+j} = \rho$ belongs to $\fc$.
Taking into account these considerations we find
that the weights contributing to the sum in \paros\ are:
$$
\eqalign{
&\nu_i =\epsilon(\omega_i) \cdot (n + 1 -i,1
\dots \buildrel i \over 2 \dots 1), \cr
& \nu_{l+1} = \epsilon(\omega_{l+1}) \cdot \rho, \cr
&\nu_{l+1+i} = \epsilon(\omega_{l+1+i}) \cdot \rho, \cr }
\qquad
\eqalign{
& i = 1, \dots, n, \cr
& \; \cr
&n = 2i-1.  \cr }
\eqn\pepa$$
where $\epsilon(\omega)$ is the
signature of the Weyl chain, given by the number
of Weyl reflections we have made to bring the weights to this form:
$$\eqalign{
&\omega_i = \sigma_1 \dots \sigma_{i-1}, \cr
& \omega_{l+1} =I, \cr
&\omega_{l+1+i} = \sigma_{l+1-i} \dots \sigma_{l-1} \sigma_l \sigma_{l-1} \dots
\sigma_{l-i+1}, \cr}
 \qquad
\eqalign{
&\Rightarrow \;\; \epsilon(\omega_i) = (-1)^{i-1}, \cr
&\Rightarrow \;\; \epsilon(\omega_{l+1}) = (-1)^0 = 1, \cr
&\Rightarrow \;\; \epsilon(\omega_{l+1+i}) = (-1)^{2i-1} = -1. \cr}
$$

 Using these results and
the scalar products in (A.22) the sum in \paros\ becomes,
$$
W_{\Lambda}^{(n,-m)}  | \rho \rangle
=\sum\limits_{i=1}^n t^{- {nm \over 2}  - m(2l+1-2i)} (-1)^{i-1}| \nu_i \rangle
+  \Biggl \{ {{0, \;\ \, \,\ n\,\ \rm odd;} \atop {
| \rho \rangle, \;\ \, \;\ n\,\ \rm even. }}
\eqn\citerea
$$
 This
equation is valid for any $n\geq 1$ as long as $l>n$.
The vacuum expectation value \vacio\ which enters \equis\ takes the form:
$$
\eqalign{
V^{(n,-m)}_{\lambda^{(1)}} & =  {\langle\rho|S W^{(n,-m)}_{\lambda
_1}|\rho\rangle \over \langle\rho |S|\rho\rangle }  \cr
& = \sum\limits_{i=1}^n t^{- {nm \over 2}  - m(2l+1-2i)} (-1)^{i-1}
 {S_{\rho ,\nu_i} \over S_{\rho , \rho }}
 + \Biggl \{ {{0, \;\ \, \,\ n\,\ \rm odd;} \atop {
 1, \;\ \, \;\ n\,\ \rm even. }} \cr}
 \eqn\isabelle
$$
The weights $\nu_i$ have the general expresion $\nu_i = \rho +
 (n -i,0 \dots \buildrel i \over 1 \dots 0)= \rho + \Lambda$.
If $\Lambda$ is a highest weight, the ratio
$S_{\rho,\rho+\Lambda}/S_{\rho,\rho}$ can be written in terms of the character
associated to $\Lambda$ with the help of the Weyl formula,
$$
{S_{\rho,\rho + \Lambda}\over S_{\rho,\rho}}={ \sum_{w \in
W} \epsilon(w) t^{ \rho \cdot w(\rho+\Lambda)} \over  \sum_{w
\in W} \epsilon (w) t^{ \rho \cdot w(\rho)} }={\rm ch}_{\Lambda}
[-{2 \pi i \over 2 k+\cox} \rho].
\eqn\cecilia
$$
All the weigths entering
\isabelle\ can be thought as highest weights and therefore we
can express $V^{(n,-m)}_{\lambda^{(1)}}$ in terms of characters:
$$
V^{(n,-m)}_{\lambda^{(1)}}=
\sum\limits_{i=1}^n t^{- {nm \over 2}  - m(2l+1-2i)} (-1)^{i-1}
   {\rm ch}_{(n-i)\lambda^{(1)} +\lambda^{(i)}} [-{2 \pi i\over 2k+ \cox}\rho]
 + \Biggl \{ { {0, \;\ \, \,\ n\,\ \rm odd;} \atop
{ 1, \;\ \, \;\ n\,\ \rm even. }}
\eqn\eloisa
$$

Let us compute
first $V^{(1,0)}_{\lambda^{(1)}}$, which is the quantity entering the
denominator in \equis. From \isabelle\ and \cecilia\ follows that one needs
to compute the character
for the fundamental representation. This calculation is done
very simply just summing over the weights of the representation:
$$
{\rm ch}_{\lambda^{(1)}}
[-{2\pi i \over 2 k+\cox}\rho] = \sum_{\mu \in M_{\lambda^{(1)}}} t^{-\mu \cdot
\rho}
=1 + \sum_{j=1}^l t^{-\mu_j \cdot\rho} +
\sum_{k=1}^l t^{-\mu_{l+1+k} \cdot\rho}
= 1 + { {t^l - t^{-l}} \over { t^{1 \over 2}
 - t^{-1 \over 2}}}.
\eqn\manzana
$$
Using this result, it turns out that
$$
V^{(1,0)}_{\lambda^{(1)}}=1 + { {\lambda - \lambda^{-1}} \over { t^{1 \over 2}
 - t^{-1 \over 2}}},
\eqn\judith
$$
which has been written entirely in terms of the variables $\lambda$ and $t$ in
\limon\ (notice that in this case $N=2l+1$). This result  agrees
with previous calculations for the unknot [\horne,\haya].

For representations different than the fundamental one, however, it is more
useful to compute the character using its expression in term of a product over
positive roots:
$$
{\rm ch}_{\Lambda}
[-{2\pi i \over 2k+\cox} \rho] =
\sum_{\mu \in M_\Lambda} t^{-\mu \cdot \rho} = \prod_{\alpha>0}{t^{{1\over
2}\alpha \cdot (\rho +  \Lambda)}-
t^{-{1\over 2}\alpha \cdot (\rho+\Lambda)}\over
t^{{1\over 2} \alpha \cdot \rho}-t^{-{1\over 2}\alpha \cdot \rho}}.
\eqn\luiza
$$
In this equation, the symbol $\alpha > 0$ indicates that the product has to be
performed over all the positive roots. For $B_l$ these are given
in the Appendix.
 Our next task is to compute the characters appearing in \eloisa\ with the
help of this formula.

In order to simplify our notation, from now on we will denote
${\rm ch}_{\Lambda} [-2\pi i\rho / (2k+\cox)]$
simply by ${\rm ch}_{\Lambda}$. Also, we
introduce  the following notation regarding  $q$-numbers and $q$-factorials:
$$
\eqalign{
[p] =& t^{p \over 2} - t^{-{p\over 2}},\cr
[p]!=&[p] [p-1]\dots [1], \;\;\; [0]! = 1.\cr}
\eqn\qunumeros$$
This allows us to write the character formula in the form:
$$
{\rm ch}_\Lambda= \prod_{\alpha>0}{[\alpha \cdot (\rho +\Lambda)]\over
[\alpha \cdot  \rho]}.
\eqn\amelia
$$

In order to compute \eloisa\ we must perform
the products in \amelia\ for weigths of
the form $(n-i)\lambda^{(1)} +\lambda^{(i)}$.
Taking into account the form of the
positive roots  listed in  (A.6), this suggests to organize the product in
\amelia\ spliting the set of positive roots in two  groups, I and II,
 depending on
whether the positive root contains
the simple root $\alpha_{(1)}$ or not. Another
thing we have to take into account
is that the metric between fundamental weigths
and simple roots of this algebra,
for the normalization chosen for the long roots,
is the following:
$$\alpha_{(i)} \cdot \lambda^{(j)} ={\rm diag} \,( 1 \dots 1,{ 1 \over 2}),
\eqn\metric
$$
due to the fact that the simple root $\alpha_{(l)}$ is shorter than the others.
 Let us carry out the
computation of the character.

The products of the positive roots with the Weyl vector are:
$$
\eqalign{
& \beta_{(j)} \cdot \rho = l-j + {1 \over 2}, \cr
& \gamma_{(j,k)} \cdot \rho = 1 + k, \cr
& \delta_{(j,k)} \cdot \rho = 2l - 2j - k.  \cr}\eqn\roo$$
and with the weigths  $\nu_i = \rho
+ (n-i,0 \dots \buildrel i \over 1 \dots 0)$, $i = 1, \dots, l-1$:

a) group I, positive roots with $\alpha_{(1)}$:
$$\eqalign{
&\beta_{(1)}\cdot \nu_i = l-1 + {1 \over 2} + n-i + 1, \cr & \, \cr
&\gamma_{(1,k)}\cdot \nu_i = 1+k+ n-i + \Biggl \{ {{1, \;\;\; i\leq k+1;}
\atop
{0, \;\;\; k \leq i-2;} } \cr & \, \cr
&\delta_{(1,k)}\cdot \nu_i = 2l - 2 -k + n-i +
\Biggl \{ {{1, \;\;\; i\leq k+1;}
 \atop {2, \;\;\; k \leq i-2;} } \cr}
\eqn\uno$$

b) group II, positive roots without $\alpha_{(1)}$:
$$\eqalign{
&\beta_{(j)}\cdot \nu_i = \cases{l-j+ {1 \over 2} + 1,  &$j \leq i;$ \cr
l-j + {1 \over 2}, &$j>i;$ \cr} \cr
& \,\ \cr
&\gamma_{(j,k)}\cdot \nu_i =  \cases{ 1+k, &$j>i \;\; {\hbox{\rm or}} \;\;
  i \geq j+k+1; $ \cr
2+k,  &$j \leq i \leq j+k; $ \cr} \cr
 & \, \cr
&\delta_{(j,k)}\cdot \nu_i =\cases{
2l-2j-k, &$j > i;$ \cr
2l-2j-k+1,  &$i \leq j \leq j+k; $ \cr
2l-2j-k+2,  &$i \geq  j+k+1.$\cr} \cr}
\eqn\dos$$

 We have these two  contributions to the characters:
$$
\eqalign{
\prod\limits_{\alpha \in {\hbox{\rm I}}}
{ [\alpha \cdot \nu_i] \over [\alpha \cdot \rho] }
& ={{[l + {1 \over 2} +n-i]} \over {[l-{1 \over 2}]}} \times
\prod\limits_{k=0}^{i-2} {{[1+k + n-i]} \over {[1+k]}}
 \times \prod\limits_{k=i-1}^{l-2} {{[2+k + n-i]} \over {[1+k]}} \times \cr
&\prod\limits_{k=0}^{i-2} {{[ 2l + n-i-k]} \over {[2l-2-k]}}
 \times \prod\limits_{k=i-1}^{l-2} {{[2l-1-k+ n-i]} \over {[2l-2-k]}} \cr
& \;\ \cr & = {1 \over { [n]\,\ [2l+n-2i+1]}} \times{{[n-i+l+{1 \over 2}]}
 \over {[l-{1 \over 2}]}}
 \times {{[2l+n-i]!} \over {[n-i]!\,\ [2l-2]!}}.\cr}
\eqn\groupuno
$$
and,
$$
\eqalign{
\prod\limits_{\alpha \in {\hbox{\rm II}}}
{ [\alpha \cdot \nu_i] \over [\alpha \cdot \rho] }
& = \prod\limits_{j=2}^i {{[l-j+1+{1 \over 2}]} \over {[l-j+{1 \over 2}]}}
 \times \prod_{j=2}^{i-1} \prod\limits_{k=0}^{i-j-1}{{[2l-2j-k+2]} \over
{[2l-2j-k]}} \times \cr
&\prod\limits_{j=2}^i \prod\limits_{k=i-j}^{l-j-1}
{{[ 2l-2j-k+1]} \over {[2l-2j-k]}}
 \times \prod\limits_{j=2}^i \prod\limits_{k=i-j}^{l-j-1 }
{{[2+k]} \over {[1+k]}} \cr
& \;\ \cr & = {{[l-{1 \over 2}]} \over { [l-i+{1 \over 2}]}}
\times{{[2l - 2i+1]} \over {[i-1]!}}
 \times {{[2l-2]!} \over { [2l-i]!}}.\cr}
\eqn\groupdos
$$
Taking into account \groupuno\  and \groupdos\ we finally
obtain a formula for the character in terms of
$q$-numbers:

$$\eqalign{
{\rm ch}_{(n-i)\lambda^{(1)} + \lambda^{(i)}}=&
\prod\limits_{\alpha>0} { [\alpha \cdot \nu_i] \over [\alpha \cdot \rho] }
= \bigg(  {1 \over { [n]}} + { 1 \over {[n+2l-2i+1]}} \bigg) \cr
&\,\,\,\,\,\,\,\,\,\,\,\,\, \times {1 \over {[n-i]!\,\ [i-1]!}}
 \prod\limits_{j=-(i-1)}^{n-i}[2l+j]. \cr}
\eqn \donimo
$$
{}From this it is
straighforward to write an expression for \eloisa\ involving only the
variables $t$ and $\lambda$. First we introduce the notation:
$$\eqalign{
[p;y] =& t^{ p \over 2}\lambda^{y} -  t^{- p \over 2}\lambda^{-y}, \cr
\beta =& n-i ,\cr
\gamma =& i-1. \cr}
\eqn\defini
$$
Recall  that  $\lambda$  is defined as $\lambda=t^{N-1\over 2}= t^l$.
One finds,
$$
\eqalign{
V_{\lambda^{(1)}}^{(n,-m)} = &
 \sum\limits_{\gamma +\beta +1=
n\atop \gamma,\beta\geq 0} t^{- {m \over 2}(\beta -
\gamma)} \lambda^{-m}
   (-1)^{\gamma}
  \bigg(
{1 \over { [n]}} + { 1 \over { [ \beta - \gamma ; 1 ]}} \bigg) \cr
&\,\,\,\,\,\,\,\,\,\,\,\,\,\,\,
\times {1 \over {[\beta]!\,\ [\gamma]!}}
 \prod\limits_{j=-\gamma }^{\beta }
[j;1] + \cases{ 0, &n odd; \cr 1, &n even. \cr}
  \cr}
\eqn\elena
$$
It remains only to obtain the deframing phase factor. The conformal weight for
the fundamental representation of $SO(2l+1)$ is given by \cyh:
$$
h_{\rho+\lambda^{(1)}}={(\rho+\lambda^{(1)})^2-\rho^2 \over 2(2k+\cox)}= {l
\over (2k+\cox)},
\eqn\manza
$$
which gives the deframing factor:
$$
\ex^{2\pi i nm h_{\rho+\lambda^{(1)}}} = t^{lmn}=\lambda^{nm}.
\eqn\pera
$$
{}From \elena\ \pera\  and \judith\
one obtains the final expression for the knot
invariant \equis, which equals the one stated in Theorem 3.1.
This ends the  proof for the case $SO(2l+1)$.

\section{$SO(2l)$}

As the calculation procedure is the same  as in the previous case,
we will simply give the main results at each step. The Lie
algebra is now $D_l$ and
 its main features are summarized in the Appendix.

The action of the knot operator on the vacuum state is given by:
$$
W_{\Lambda}^{(n,-m)}  | \rho \rangle =
\sum_{i=1}^{2l} t^{- {1 \over 2}
\mu_i^2 nm - m \mu_i \cdot \rho} | \rho + n \mu_i \rangle ,
\eqn\ramboia
$$
where $\mu_i = 1,...,2l$
are the weigths in $M_{\lambda^{(1)}}$ whose expression is in (A.25).
The vectors $\rho +n\mu_i$ have the structure:
$$
\eqalign{
&\rho + n\mu_1 = (n+1,1,\dots,1), \cr
&\;\;\;\;\;\;\; \vdots \cr
&\rho + n\mu_j =  (1,\dots,1,1-n,\buildrel j \over {1+n},1,\dots,1),\cr
&\;\;\;\;\;\;\;  \vdots \cr
&\rho + n\mu_{l-1} = (1,\dots,1,1-n,1+n,1+n), \cr
&\rho + n\mu_l = (1,\dots,1,1-n,1+n), \cr
 & \rho + n\mu_{l+1} = (1,\dots,1,1+n,1-n), \cr
& \rho + n\mu_{l+2} = (1,\dots,1,1+n,1-n,1-n), \cr
& \;\;\;\;\;\;\;  \vdots \cr
& \rho + n\mu_{l+j} =  (1,\dots,1,1+n,\buildrel l-j \over {1-n},1,\dots,1),
 \cr
 & \;\;\;\;\;\; \vdots \cr
& \rho + n\mu_{2l-1} = (1+n,1-n,1,\dots,1), \cr
&  \rho + n\mu_{2l} = (1-n,1,\dots,1), \cr}
\eqn\pesitos$$
and those who contribute
to the sum in \ramboia, for the case $n<l$, after taking the  suitable
chain of Weyl reflections, and the corresponding signature, are:
$$
\eqalign{
&\nu_i =(-1)^{(i-1)} \cdot
(n + 1 -i,1, \dots, \buildrel i \over 2, \dots, 1), \cr
&\nu_{l+1+i} = \rho, \cr }
\qquad
\eqalign{
& i = 1, \dots, n, \cr
&n = 2i.  \cr }
\eqn\contribuyen
$$
We see that, very similarly to the
$SO(2l+1)$ case, the number of weigths we have to
take into account is bounded
by $n$ and that there is an extra one in the case of
$n$ even.
So the expression \ramboia\
becomes, after using the scalar products in (A.27):
$$
W_{\Lambda}^{(n,-m)}  | \rho \rangle =
 \sum\limits_{i=1}^n t^{- {nm \over 2}  - m(l-i)} (-1)^{i-1}| \nu_i \rangle
 + \Biggl \{ {{0\, ,\;\;\; n\; \rm odd;} \atop {| \rho \rangle\, , \;\,  n \;
\rm even. }}
\eqn\argos
$$
The quantity
$V^{(1,0)}_{\lambda^{(1)}}$ is obtained from the character of the fundamental
representation:
$$
V^{(1,0)}_{\lambda^{(1)}} = \sum_{\mu \in M_{\lambda^{(1)}}} t^{-\mu \cdot
\rho} = 1 +
{{ \lambda - \lambda^{-1} } \over {t^{ 1 \over 2} - t^{-{1 \over 2}}}}=
{{ [1] + [0;1]} \over [1]},
\eqn\nonudo$$
where for the last equality
we have use the definitions \qunumeros\ and \defini.
To calculate $V_{\lambda^{(1)}}^{(n,-m)} $
we again need the characterization of the positive roots which is
contained in (A.7).
Then one computes the products of these roots with
the Weyl vector and the weigths $\nu_i$. From these we obtain the following
formula for the characters:
$$
\prod\limits_{\alpha>0} { [\alpha \cdot \nu_i] \over [\alpha \cdot \rho] }
= \bigg(  {1 \over { [n]}} + { 1 \over {[n+2l-2i]}} \bigg)
 {1 \over {[n-i]!\,\ [i-1]!}}
 \prod\limits_{j=-(i-1)}^{n-i}[2l+j-1].
\eqn\laertes
$$
Using again \defini, one gets:
$$
\eqalign{
V_{\lambda^{(1)}}^{(n,-m)} = &
 \sum\limits_{\gamma +\beta +1=n\atop\gamma,\beta\geq 0} t^{- {m \over 2}(\beta
- \gamma)} \lambda^{-m}
   (-1)^{\gamma}
  \bigg(  {1 \over { [n]}}
+ { 1 \over { [ \beta - \gamma ; 1 ]}} \bigg) \cr
&\,\,\,\,\,\,\,\,\,\,\,\times {1 \over {[\beta]!\,\ [\gamma]!}}
 \prod\limits_{j=-\gamma }^{\beta } [ j; 1] +
\cases{ 0, &n odd; \cr 1, &n even. \cr}
 \cr}
\eqn\naxos
$$
The framing factor for this case is given by:
$$
\ex^{2 \pi i nm h_{\rho +
\lambda^{(1)}}} = t^{{{nm} \over 2} (2l - 1)} = \lambda^{nm}.
\eqn\prefactor
$$
It is easy to see that
taking into account \prefactor, \nonudo, and \naxos\ we obtain
the expression \viola\
for the knot invariant associated to the fundamental representation
of $SO(2l)$. Notice also that although $\lambda$ is defined in a different way
with respect to the
rank of the algebra, $l$, its definition is the same for both cases
in terms of  the variable
$N$ of $SO(N)$. This completes the proof of  Theorem 3.1.

\section{Natural variables of the Kauffman polynomial and Yokota's formula}

The Dubrovnik version of the Kauffman polynomial, as described in pag. 215
of \REF\kyp{L.H. Kauffman, ``Knot and Physics", World Scientific, 1991}
[\kyp], depends
on two variables, $a$ (which is called $\alpha$ in [\kyp])
and $z$. We will  refer to these variables as the natural ones.
In those variables the skein rules have the simple form
shown in [\kyp].
We will denote  the Dubrovnik version of the Kauffman polynomial, normalized in
such a way that for the unknot its value is one, by $Y_K(a,z)$,  and will try
to
identify these  variables in terms of ours. This can be done comparing the
skein rules in [\kyp]  to the skein rules obtained from Chern--Simons
theory in  [\ygw,\horne,\kcp]. It turns out that,
$$
\eqalign{
a=&\lambda = e^{2 \pi i  h_{\rho + \lambda^{(1)}}},\cr
z=&[1]=t^{1\over 2}-t^{-{1\over 2}}.}
\eqn\naturales
$$
The formula in Theorem 3.1 can therefore be stated as:
$$
Y_{n,m}(\lambda,t^{1\over 2} - t^{-{1\over 2}}) =
X_{\lambda^{(1)}}^{(n,m)},
\eqn\reme
$$
where $X_{\lambda^{(1)}}^{(n,m)}$ is given in \viola.

To compare our formula \reme\ to the one obtained by Yokota in
[\yoko] we will use \naturales\ and the identification done
in [\yoko] between its variables, $q$ and $\alpha$, and the natural ones.
Proceeding in this way one concludes that the relation between
our variables and Yokota's is:
$$
\eqalign{
q &= t^{-{1 \over 2}}, \cr
\alpha^2 &= -(q\lambda)^{-1}. \cr}
\eqn\yokotas
$$
Taking into account that Yokota uses an orientation opposite to
ours, and therefore we must compare \reme\ to its formula for
$Y_{n,m}(a^{-1},-z)$, one finds complete agreement after substituting
\yokotas\ in the formula given in [\yoko].

\endpage

\chapter{Relation between the HOMFLY  and Kauffman polynomials for torus knots}

The HOMFLY [\homfly]
  and Kauffman polynomials [\kau] have the common characteristic of being
functions of two variables defined for oriented links, although their behavoir
under change of orientation of some of the link components is quite different.
On the other hand, the skein rules that define them are also different:
in the first one the relation is established among three diagrams
and in the second one among four.
Both are able to  differenciate in many cases one knot from its mirror image,
although  Kauffman's is more powerfull in this sense. These two polynomials are
considered as independent, in the sense that there is not a subtle
change of variables taking one into the other. In
\REF\nudos{W.B.R. Lickorish and K.C.Millett \journal\mm &61(88)}[\nudos]
there are examples of knots with the same Kauffman and different
HOMFLY and {\it viceversa}. We will prove that for the
 particular case of torus
knots there is a relation between these two polynomials.
Let's begin recalling the expression of the HOMFLY polynomial for torus
knots. It was first obtained in
[\jonesAM], reobtained in
[\rosso] using quantum groups and in
[\lama] from the Chern-Simons theory with gauge group $SU(N)$.
 The corresponding invariant has the form [\lama]:
$$
\eqalign{
P_{n,m} (a,z) &= P_{n,m}
((\lambda t)^{1\over 2},t^{1\over 2}- t^{-{1\over 2}}) \cr &
= \Big ( {1-t \over 1-t^n }\Big)
{\lambda^{{1\over2}(m-1)(n-1)} \over \lambda t-1} \sum_{p+i+1=n \atop p, i
\ge 0} (-1)^i t^{mi + {1 \over 2} p(p+1)} {\prod_{j=-p}^i (\lambda t-t^j)
\over \prod_{j=1}^i (t^j-1) \prod_{j=1}^p (t^j-1)}\cr
&= {[1] (\lambda t)^{{1\over 2}m(n-1)} \over [-1;-{1\over 2}]}
\sum_{\beta+\gamma+1=n\atop \beta,\gamma \geq 0}(-1)^\gamma t^{{m\over
2}(\beta-\gamma)} {1 \over [n] [\beta]!
[\gamma]!}\times \prod_{j=-\gamma}^{\beta} [j-1;-{1\over 2}],\cr}
\eqn\nuria
$$
where:
$$
\eqalign{
\lambda &= t^{N-1},\cr
t & = \ex^{2\pi i \over k+ g^{ \vee}}.\cr}
\eqn\gatito
$$
If one performs one of these two changes of variables :
$$
 t^{1 \over 2} \rightarrow t^{-{1 \over 2}},
\,\,\,\,\,\,\,\,\,\,\,  {\hbox{\rm or}}  \,\,\,\,\,\,\,\,\,\,\,\,
 t^{1 \over 2} \rightarrow -t^{1 \over 2},
\eqn\changeses
$$
one finds that \viola\ transforms into:
$$
\eqalign{
Y_{n,m}(a,-z)&
=Y_{n,m}(\lambda,t^{-{1\over 2}}-t^{1\over 2}) = \cr& =-{{[1] \lambda^{nm}}
\over { {[1] - [0;1] }  }} \times \Biggl(
 \sum\limits_{\gamma +\beta +1=n \atop \gamma,\beta\geq 0} t^{- {m \over 2}
(\beta - \gamma)} \lambda^{-m} (-1)^{\gamma}
 \times \bigg(  {1 \over { [n]}} -
{ 1 \over { [ \beta - \gamma ; 1 ]}} \bigg) \cr
\,\,\,\,\,\,\,\,\,\,\,\,\,\,
&\,\,\,\,\,\,\,\,\,\,\,\,
\times {1 \over {[\beta]!\,\ [\gamma]!}} \times
 \prod\limits_{j=-\gamma }^{\beta } [j; 1] +
\cases{ 0, &n odd; \cr -1, &n even. \cr}
\;\ \Biggr), \cr}
\eqn\reofelia
$$
It is worth to remark that this is exactly the formula obtained
when one calculates
the polynomial for torus knots associated to the fundamental representation
of $Sp(N)$ from Chern--Simons theory. This can be shown explicitly
using the methods developed in the previous section, or
from the form of the skein rules for the fundamental of
$Sp(N)$ obtained in [\horne,\kcp] from Chern-Simons gauge theory.
Let us compare \viola, \nuria, and \reofelia. The crucial point is that, using
the auxiliary variable $q$, these
three expressions can be written as follows:
$$
\eqalign{
Y_{n,m}&(a,q-q^{-1}) = \cr &={ a^{nm} [1]_q \over
[1]_q + a-a^{-1}} \times \Biggl(
\sum\limits_{\gamma +\beta +1=n\atop \gamma,\beta\geq 0}
q^{- {m}(\beta - \gamma)} a^{-m} (-1)^{\gamma}
 \times \bigg(  {1 \over { [n]_q}} + { 1 \over {
q^{\beta-\gamma}a-q^{\gamma-\beta}a^{-1}}} \bigg) \cr
&\,\,\,\,\,\,\,\,\,\,\,\,\,\,\,\,\,\, \times {1 \over {[\beta]_q!\,\
[\gamma]_q!}} \times \prod\limits_{j=-\gamma }^{\beta }
(q^j a - q^{-1} a^{-1}) + \cases{ 0, &n odd; \cr 1, &n even; \cr}
 \Biggr), \cr}
\eqn\naxosdos$$
$$
\eqalign{
Y_{n,m}&(a,-(q^{-1}-q))  = \cr &= -{ a^{nm} [1]_q \over
[1]_q - a+a^{-1} } \times \Biggl(\sum\limits_{\gamma +\beta +1=n \atop
\gamma,\beta\geq 0} q^{- {m}(\beta - \gamma)} a^{-m} (-1)^{\gamma}
\times \bigg(  {1 \over { [n]_q}} -
{ 1 \over { q^{\beta-\gamma}a-q^{\gamma-\beta}a^{-1}}} \bigg) \cr
&\,\,\,\,\,\,\,\,\,\,\,\,\,\,\,\,\,\,
\times {1 \over {[{\beta}]_q!\,\ [{\gamma}]_q!}}
 \times \prod\limits_{j=-\gamma }^{\beta }
 (q^ja-q^{-j}a^{-1}) +
\cases{ 0, &n odd; \cr -1, &n even; \cr}
\;\ \Biggr), \cr}
\eqn\ofeliados
$$
and,
$$
\eqalign{
P_{n,m}&(a,q^{-1}-q)  = \cr &= { a^{m(n-1)} [1]_q \over
a-a^{-1} } \sum\limits_{\gamma +\beta +1=n \atop
\gamma,\beta\geq 0} q^{- {m}(\beta - \gamma)}  (-1)^{\gamma}
{1 \over [n]_q \, {[{\beta}]_q!\,\ [{\gamma}]_q!}}
\prod\limits_{j=-\gamma }^{\beta }
 (q^ja-q^{-j}a^{-1}), \cr}
\eqn\nuriados
$$
where,
$$
[n]_q = q^n-q^{-n}.
\eqn\nuevoq
$$
The structure on the right hand side of \nuriados\  shows  that the
 HOMFLY polinomial, $P_{n,m}(a,z)$,  can be expresed in terms of a linear
 combination of the polynomials $Y_{n,m}(a,z)$ and $Y_{n,m}(a,-z)$.
In fact, after performing some algebra from
 \naxosdos, \ofeliados\ and \nuriados, one obtains,
$$
P_{n,m}(a,z)={1\over 2}(Y_{n,m}(a,z)+Y_{n,m}(a,-z)) +
{z\over 2(a-a^{-1})}(Y_{n,m}(a,z)-Y_{n,m}(a,-z)).
\eqn\gaussdosdos
$$
This ends the proof of the relation \gaussdos\ between the
 HOMFLY and Kauffman polynomials which was presented in the
introduction.

For $a=1$, the ordinary version of the Kauffman polynomial, $F_K(a,z)$, becomes
the unoriented polynomial invariant of ambient isotopy discovered in
\REF\blm{R.D. Brandt, W.B.R. Lickorish and K.C. Millet,
{\sl Invent. Math.} {\bf 84} (1986) 563}
\REF\ho{C.F. Ho, AMS Abstract, vol. 6, no. 4, Issue {\bf 39}
(1985) 300} [\blm,\ho] which is usually denoted by $Q_K(z)=F_K(1,z)$.
Similarly, we define:
$$
\tilde Q_K(z)=Y_K(1,z).
\eqn\mcarmen
$$
It turns out that for torus knots, after performing the limit
$N\rightarrow 1$ in \viola, which is equivalent to $a\rightarrow 1$,
one finds:
$$
\tilde Q_{n,m}(z) = 1.
\eqn\mjesus
$$
In the case of the HOMFLY polynomial, the limit $a\rightarrow 1$
leads to the Alexander-Conway polynomial, $\Delta_K(z)=P_K(1,z)$.
{}From \gaussdosdos\ and \mjesus\ one finds:
$$
\Delta_{n,m}(z) = 1 + {z\over 4}{\partial\over\partial a}
\Big(Y_{n,m}(a,z)-Y_{n,m}(a,-z)\Big)\Big|_{a=1}.
\eqn\olga
$$
Notice that this expression is consistent with the fact that
$\Delta_{n,m}(z)$ must be 1 plus a polynomial containing only even
powers of $z$.

\endpage
\chapter{Conclusions and prospects}

In this paper we have presented the construction of the
operator formalism, originally discussed in \nos\ and  \llr\
for the groups $SU(2)$ and $SU(N)$ respectively, for
an arbitrary simple group. The main result in this respect
is the general form for knot operators presented in
\venator.

Knot operators are utilized to compute the knot invariant corresponding
to the fundamental representation of the gauge group
$SO(N)$. The resulting formula is presented in \reme\
and \viola, and shown to agree with a previous expression
for the Kauffman polynomial. This formula is compared to
known expressions for the HOMFLY polynomial and the  relation
\gaussdos\
between the Kauffman and the HOMFLY polynomials for torus knots
is proved.

Our result \gaussdos\ confirms that the Kauffman polynomial is
more fundamental than
the HOMFLY polynomial. The simplicity of the relation
obtained suggest that
it could be obtained by other methods. In this respect it would
very interesting if it could be reobtained using skein rules.

It would be also worthwhile to study how our results fit in Jaeger's
expansions for
 the Kauffman polynomial in terms of HOMFLY polynomials (see
for example pag. 219 of
[\kyp])
Finally, one should also study if there exist similar formulas for other sets
of knots. In this respect one would like to start studying the situation in
sets characterized by a generalization of the notion of a torus knot.
A torus knot is a knot that can be placed on the surface of a
standardly embedded torus
in $S^3$ without self-intersection. There are knots which
can be placed on a standardly
embedded genus two surface without self-intersection
but not on a genus one surface. One could analyze for example if there is
a relation of the type \gaussdos\
for these knots. In general one could study the
problem for knots placed on a genus $g$ surface. Work in
this direction will be presented elsewhere.

\vskip0.5cm

\ack

We would like to thank A. V. Ramallo and A. Saa for very helpful
discussions and comments. We would also like to thank M. Mari\~no
for a careful reading of the manuscript.
This work was supported in part by DGICYT under grant
PB93-0344 and  by CICYT under grant AEN94-0928.

\endpage

\appendix

\noindent{\caps Group-theoretical Conventions.}

In this section of the  Appendix we will summarize our
group-theoretical  conventions. Let $G$  be a
compact simple group  of rank $l$,
with  generators $T^a$, $a=1,\dots,{\hbox{\rm dim}}(G)$, which are chosen
to be antihermitian. For the fundamental representation of $G$ they
are normalized  as follows:
$$
\tr(T^a T^b)=-y\psi^2\delta^{ab}
\eqn\ai
$$
where $y$ is the Dynkin index of the fundamental representation
and $\psi^2$ is the squared length
of the longest simple root of $G$.
The value of $y$ for the  groups $SU(N)$, $SO(N)$, $Sp(N)$,
$E_6$, $E_7$, $E_8$, $F_4$ and $G_2$ are $1/2, 1, $1/2, 9, 12, 30,
6 and 3, respectively.

We will denote the  $l$ fundamental roots of $G$  by
$\alpha_i$, $i=1,...,l$. In the explicit calculations carried out
in sect. 3 they have been chosen
in such a way that the long roots
have length $\sqrt{2}$, \ie, $\psi^2=2$.
The  Cartan matrix $g_{ij}$,
$$
g_{ij}=2 {\alpha_{(i)}\cdot\alpha_{(j)} \over
\alpha_{(i)}\cdot\alpha_{(i)}},
\eqn\cartan
$$
takes the following forms for the two Lie algebras
$B_l$ ($l=(N-1)/2$, $N$ odd) and $D_l$ ($l=N/2$, $N$ even)
associated to the simple group
$SO(N)$, which is the one that has been considered in
this paper:
$$
g_{ij}(B_l)=\left(\matrix{2&-1&0&0&\cdot&\cdot&\cdot&\cdot&0\cr
                -1&2&-1&0&\cdot&\cdot&\cdot&\cdot&0\cr
                0&-1&2&-1&\cdot&\cdot&\cdot&\cdot&0\cr
                \cdot&\cdot&\cdot&\cdot&\cdot&\cdot&\cdot&\cdot&\cdot\cr
                \cdot&\cdot&\cdot&\cdot&\cdot&\cdot&\cdot&\cdot&\cdot\cr
                0&0&0&0&\cdot&\cdot&2&-1&0\cr
                0&0&0&0&\cdot&\cdot&-1&2&-1\cr
                0&0&0&0&\cdot&\cdot&0&-2&2\cr}\right),
\eqn\aii
$$
and,
$$
g_{ij}(D_l)=\left(\matrix{2&-1&0&0&\cdot&\cdot&\cdot&\cdot&\cdot&0\cr
                -1&2&-1&0&\cdot&\cdot&\cdot&\cdot&\cdot&0\cr
                0&-1&2&-1&\cdot&\cdot&\cdot&\cdot&\cdot&0\cr
                \cdot&\cdot&\cdot&\cdot&\cdot&\cdot&\cdot&\cdot&\cdot&\cdot\cr
                \cdot&\cdot&\cdot&\cdot&\cdot&\cdot&\cdot&\cdot&\cdot&\cdot\cr
                0&0&0&0&\cdot&\cdot&\cdot&\cdot&0&0\cr
                0&0&0&0&\cdot&\cdot&-1&2&-1&-1\cr
                0&0&0&0&\cdot&\cdot&0&-1&2&0\cr
                0&0&0&0&\cdot&\cdot&0&-1&0&2\cr}\right).
\eqn\aij
$$

We will denote the root lattice by $\rl$.
This $l$-dimensional space is generated by
the fundamental roots $\alpha_{(i)}$, which
can be taken as a basis, the root basis.
Any vector $x$  in this basis has components
${x}^i$ given by:
$$
x=\sum_{i=1}^{l} {x}^i\alpha_{(i)}.
\eqn\aiip
$$
Among all the roots in $\rl$ there is a subset which plays an important role in
the calculation performed in the paper. These are the positive roots.
For $SO(N)$ they take the form
\REF\corn{J. F. Cornwell, {\it Group theory in physics}, Vol. 3,
Academic Press, 1989 }[\corn],

\noindent -  algebra $B_l$:
$$\eqalign{
&\beta_{(j)} = \alpha_{(j)} + \dots  + \alpha_{(l)},
\,\,\,\,\,\,\,\,\;\;\;\;\;\;\;\; j = 1, \dots, l, \cr
& \gamma_{(j,k)} = \alpha_{(j)} +\dots + \alpha_{(j+k)},
\;\;\;\;\;\;\;\; i = 1, \dots, l-1, \;\;\;\ k = 1, \dots, l-j-1, \cr
& \delta_{(j,k)} = \alpha_{(j)} + \dots + \alpha_{(j+k)} +
2( \alpha_{(j+k+1)} + \dots + \alpha_{(l)}), \cr
& \;\;\;\;\;\;\;\;\;\;\;\;\;\;\;\;\;\;\;\;
\;\;\;\;\;\;\;\;\;\;\;\;\;\;\;\;\;\;\;\;\;\;\;
\;\;\;\;\;\;\;\;\;\;\;\;\;\;\;\;\;\;\;\;\; j = 1,\dots, l-1,
 \;\;\; k = 0, \dots, l-j-1. \cr}
\eqn\liebre
$$

\noindent -  algebra $D_l$:
$$
\eqalign{
& \alpha_{(j)}, \,\,\,\,\,\,\,\,\,\,\,
\;\;\;\;\;\;\;\;\;\; j = 1,\dots,l, \cr
&\beta_{(j)} = \alpha_{(j)} + \dots + \alpha_{(l-2)} + \alpha_{(l)},
\;\;\;\;\;\;\;\; j = 1, \dots, l-2, \cr
& \gamma_{(j,k)} = \alpha_{(j)} +\dots + \alpha_{(j+k)},
\;\;\;\;\;\;\;\; j = 1, \dots, l-2, \;\;\;\ k = 1, \dots, l-j, \cr
&\delta_{(j,k)} = \alpha_{(j)} + \dots + \alpha_{(j+k)} +
2 \big( \alpha_{(j+k+1)} + \dots + \alpha_{(l-2)} \big) +
 \alpha_{(l-1)} + \alpha_{(l)}, \cr
&\;\;\;\;\;\;\;\;\;\;\;\;\;\;\;\;\;\;\;\;\;\;\;\;\;\;\;\;
\;\;\;\;\;\;\;\;\;\;\;\;\;\;\;\;\;\;\;\;\;\;\;\;\;\;\;\;\;
\;\;\;\;\;\;\;\;\; j = 1,\dots, l-3, \;\;\;\ k = 0, \dots, l-3-j. \cr}
\eqn\gato
$$

The fundamental weights $\lambda^{(i)},$ $i=1,...,l,$ satisfy:
$$
2\,{\alpha_{(i)}\cdot\lambda^{(j)}
\over{\alpha_{(i)}\cdot\alpha_{(i)}}}=\delta_{i}^{j}.
\eqn\aiv
$$

 The fundamental weights generate over ${\bf Z}$ an $l$-dimensional
lattice called the weight lattice which will be denoted by $\wl$.
The lattices $\rl$ and $\wl$ are dual to each other and
$\rl\in\wl$.
The $l$-dimensional basis expanded by the fundamental weights is
called the Dynkin basis. Any vector $x$ has in this
basis components $x_i$ given by:
$$
x=\sum_{i=1}^{l}  x_i\lambda^{(i)}.
\eqn\avi
$$
The matrix $G^{ij}=\lambda^{(i)}\cdot\lambda^{(j)}$
gives the metric in weigth space,
so it allows us to rise indices.
Its expression for the  algebras $D_l$ and $B_l$  is:
$$
G^{ij}(D_l)= {1 \over 2} \left(\matrix{2&2&2&\cdot&\cdot&\cdot&2&1&1\cr
                2&4&4&\cdot&\cdot&\cdot&4&2&2\cr
                2&4&6&\cdot&\cdot&\cdot&6&3&3\cr
                \cdot&\cdot&\cdot&\cdot&\cdot&\cdot&\cdot&\cdot&\cdot\cr
                \cdot&\cdot&\cdot&\cdot&\cdot&\cdot&\cdot&\cdot&\cdot\cr
                2&4&6&\cdot&\cdot&\cdot&2(l-2)&l-2&l-2\cr
                1&2&3&\cdot&\cdot&\cdot&l-2&l/2&(l-2)/2\cr
                1&2&3&\cdot&\cdot&\cdot&l-2&(l-2)/2&l/2\cr}\right),
\eqn\gdl
$$
and,
$$
G^{ij}(B_l)= {1 \over 2} \left(\matrix{2&2&2&\cdot&\cdot&\cdot&2&1\cr
                2&4&4&\cdot&\cdot&\cdot&4&2\cr
                2&4&6&\cdot&\cdot&\cdot&6&3\cr
                \cdot&\cdot&\cdot&\cdot&\cdot&\cdot&\cdot&\cdot\cr
                \cdot&\cdot&\cdot&\cdot&\cdot&\cdot&\cdot&\cdot\cr
                2&4&6&\cdot&\cdot&\cdot&2(l-1)&l-1\cr
                1&2&3&\cdot&\cdot&\cdot&l-1&l/2\cr }\right).
\eqn\gbl
$$

Among the weights in $\wl$ there is one which plays an important role in
Chern-Simons theory because it can be regarded as the vacuum. This weight
is denoted by $\rho$ and all its components are one:
$$
\rho = \sum_{i=1}^l \lambda^{(i)}.
\eqn\elrho
$$

The irreducible representations of $G$ are characterized by
highest weights  $\Lambda$. Highest weights can
be written uniquely as a linear combination of fundamental weights
with non-negative integer coefficients $h_i$,
$$
\Lambda = \sum_{i=1}^l  h_i\lambda^{(i)}.
\eqn\aix
$$
The set of weights of an irreducible representation of highest weight $\Lambda$
will be denoted as $M_\Lambda$. To build this set one may use   the following
rule:
 if a weight
$\mu\in M_\Lambda$  has the $k^{th}$ Dynkin component
greater than zero (\ie, ${\mu}_k>0$), then the vectors obtained
by subtracting  $t\alpha_k$ ($t=1,...,{\mu}_k$) from $\mu$
are also elements of $M_\Lambda$.
One can start applying this rule to $\Lambda$
and then to the successive
weights obtained to build the different elements of $M_\Lambda$.
The multiplicities of each
weight can be obtained using Freudenthal's formula
\REF\book{N. Jacobson, {\it Lie Algebras,} Wiley-Interscience, New York,
1962; J. E. Humphreys, {\it Introduction to Lie Algebras and Representation
Theory}, Springer, New York, 1972.} [\book].

The Weyl group is generated by $r$ reflections $\sigma_i$, $i=1,...,l$, on
weight space
$$
x\in \wl, \,\,\,\,\,\, \sigma_i(x) = x -
{2\over \alpha_{(i)} \cdot \alpha_{(i)}}
\alpha_{(i)}(\alpha_{(i)}\cdot x).
\eqn\apnu
$$
It divides the weight lattice $\wl$ into
domains. The fundamental domain or Weyl chamber is chosen to be the one
containing all the weights $x\in\wl$ such that,
$$
\alpha_{(i)}\cdot x \geq 0.
\eqn\axvi
$$

The Weyl  character for an irreducible representation of highest weight
$\Lambda$ is defined as,
$$
\ch_\Lambda(a)=\sum_{\mu\in M_\Lambda} \ex^{a\cdot\mu},
\eqn\apuno
$$
where $a=a_i\lambda^{(i)}$. The Weyl character satisfies the equation
[\book],
$$
\ch_\Lambda(a)={\sum_{w\in W}\epsilon(w) \ex^{w(\Lambda+\rho)\cdot a}
\over \sum_{w\in W}\epsilon(w) \ex^{w(\rho)\cdot a} },
\eqn\apdos
$$
known as the Weyl character formula. When $a=-\rho$,
we have an expression for the character
[\kac] which is particularly useful:
$$
\sum_{\mu \in M_\Lambda} \ex^{-\mu \cdot \rho} =
\prod_{\alpha>0}{\ex^{{1\over 2}\alpha\cdot(\rho +
\Lambda)}-\ex^{-{1\over 2}\alpha \cdot(\rho+\Lambda)}\over \ex^{{1\over 2}
\alpha\cdot\rho}-\ex^{-{1\over 2}\alpha\cdot\rho}},
\eqn\lechuguino
$$
where $\alpha > 0$ denotes a sum over all positive roots.

An important set of weights used in this work is the one made by
Weyl-antisymmetric combinations of weights in $\wl / s\rl$ where $s$
is an arbitrary non-negative integer. This set of weights builds the
fundamental chamber ${\cal F}_s$.

\noindent{\caps Fundamental representation of $SO(2l+1)$}

 The fundamental representation of
$B_l$ is
 associated to the highest weight $\Lambda =
\lambda^{(1)} = (1,0,\dots,0)$, and the corresponding weight space is: $$
M_{\lambda^{(1)}}=\{ \mu_i: 1\le i \le 2l+1 \},
\eqn\rosi
$$
where:
$$
\eqalign{
&\mu_1=\lambda^{(1)}=(1,0,\dots,0),\cr
&\mu_j=\mu_{j-1}-\alpha_{(j-1)} = (0,\dots,-\buildrel j \over
1,1,0,\dots,0), \,\,\,\,\,\,\,\,\,\,\;\; j = 1 \dots l-1, \cr
&\mu_l=
\mu_{l-1}-\alpha_{(l-1)} = (0,\dots,0,-1,2),\cr
 &\mu_{l+1} =\mu_l-\alpha_{(l)} = 0, \cr
&\mu_{l+1+i} = -\mu_{l+1-i}, \,\,\,\,\,\,\;\;\;\;\;\;\;\; i = 1, \dots,
 l. \cr}\eqn\butano
$$
  We can write these weights as follows:
$$
\eqalign{
&\mu_j = \sum\limits_{i=1}^l \bigl[ -\delta_{j-1,i} + \delta_{j,i} \bigr]
 \lambda^{(i)}, \;\;\;\;\;\; j = 1, \dots, l, \;\;\;\;\;\; j \not= l-1, \cr
&\mu_{l-1} = \bigl[ -\delta_{l-1,i} + 2\delta_{l,i}
 \bigr]\lambda^{(i)}, \cr
&\mu_{l+1+i} = -\mu_{l+1-i}, \;\;\;\;\;\; i = 1, \dots, l. \cr}
\eqn\abuela
$$

 We also need the scalar products $\rho\cdot\mu_i$ and
$\mu_i\cdot\mu_i$. Using the form \butano\ and \gbl, we can easily find:
$$
\eqalign{
&\mu_i^2 = 1, \;\;\;\;\;\;\; i = 1, \dots, 2l+1, \;\;\;\;\;\; i\not= l+1, \cr
&\mu_{l+1}^2 = 0, \cr
&\rho \cdot \mu_i ={1 \over 2}[2l-(2i-1)], \;\;\;\;\;\;\; i = 1,
\dots, l, \cr
&\rho \cdot\mu_{l+1} = 0, \cr
&\rho \cdot \mu_{l+1+i} = -{1 \over 2}(2i-1), \;\;\;\;\;\;\; i = 1, \dots,
 l. \cr}
\eqn\vestido
$$
The action of the Weyl reflections on the fundamental weights
$\lambda^{(i)}$ follows from \apnu:
$$
\eqalign{
\sigma_1(x)=&(-x_1,x_2 + x_1,x_3, \dots,x_l), \cr
\sigma_i(x)=&(x_1,\dots,x_{i-1} + x_i,-x_i,x_i + x_{i+1}, \dots,x_l),
\,\,\,\,\,\,\,\;\;\; i = 1, \dots, l-2, \cr
\sigma_{l-1}(x)=&(x_1,\dots,x_{l-3},x_{l-2} + x_{l-1} , -x_{l-1}
,x_l + 2 x_{l-1}), \cr
\sigma_l(x)=&(x_1,\dots,x_{l-2},x_{l-1} + x_l , -x_l ). \cr}
\eqn\naranja
$$

\noindent{\caps Fundamental representation of $SO(2l)$}

In this section we present the results
concerning the fundamental representation of
$D_l$. It is
 associated to the highest weight $\Lambda =
\lambda^{(1)} = (1,0,\dots,0)$, and the corresponding weight space is: $$
M_{\lambda^{(1)}}=\{ \mu_i: 1\le i \le 2l \},
\eqn\paloma
 $$
where:
$$
\eqalign{
&\mu_1=\lambda^{(1)}=(1,0,\dots,0),\cr
&\mu_j= \mu_{j-1}-\alpha_{(j-1)} = (0,\dots,-\buildrel j \over
1,1,0,\dots,0),\;\;\,\,\,\,\,\,  j = 1, \dots, l-2, \cr
&\mu_{l-1}=
\mu_{l-2}-\alpha_{(l-2)} = (0,\dots,0,-1,1,1),\cr
& \mu_l =\mu_{l-1}-\alpha_{(l-1)} =
(0,\dots,0,-1,1),\cr
&\mu_{l+i} = -\mu_{l+1-i}, \;\;\;\;\;\;\; i = 1, \dots, l.\cr}
\eqn\leila
$$
  We can write these weights as follows:
$$
\eqalign{
&\mu_j = \sum\limits_{i=1}^l \bigl[ -\delta_{j-1,i} + \delta_{j,i} \bigr]
 \lambda^{(i)}, \;\;\; \;\;\; j = 1, \dots, l \,\,\,\;\;\;\; j \not= l-1, \cr
&\mu_{l-1} = \bigl[ -\delta_{l-2,i} + \delta_{l-1,i} + \delta_{l,i}
 \bigr]\lambda^{(i)}, \cr
&\mu_{l+i} = -\mu_{l+1-i},  \;\;\;\;\;\; i = 1, \dots, l. \cr}
\eqn\mortadelo
$$

 We also need the scalar products $\rho\cdot\mu_i$ and
$\mu_i\cdot\mu_i$. Using the form \mortadelo\ and \gdl, we easily find:
$$
\eqalign{
&\mu_i^2 = 1, \;\;\;\; i = 1, \dots, 2l,\cr
&\rho \cdot \mu_i = l-i, \;\;\;\; i = 1, \dots, l, \cr
&\rho \cdot \mu_{l+i} = -(i-1), \;\;\;\; i = 1, \dots, l. \cr}
 \eqn\clarisse
$$
The action of the Weyl reflections on the fundamental weights
$\lambda_i$ follows from \apnu:
$$
\eqalign{
\sigma_1(x)=&(-x_1,x_2 + x_1,x_3, \dots,x_l), \cr
\sigma_i(x)=&(x_1,\dots,x_{i-1} + x_i,-x_i,x_i + x_{i+1}, \dots,x_l)
\;\;\;\,\,\,\,\, i = 1, \dots, l-3, \cr
\sigma_{l-2}(x)=&(x_1,\dots,x_{l-4},x_{l-3} + x_{l-2} , -x_{l-2}
,x_{l-1} + x_{l-2},x_l + x_{l-2}), \cr
\sigma_{l-1}(x)=&(x_1,\dots,x_{l-3},x_{l-2} + x_{l-1} , -x_{l-1}
,x_l ), \cr
\sigma_l(x)=&(x_1,\dots,x_{l-3},x_{l-2} + x_l , x_{l-1}
,-x_l ). \cr}
\eqn\tatiana
$$

\noindent{\caps Theta functions of level $s$}

The Theta functions of level $s$ (being $s$ an arbitrary positive integer)
play a fundamental role in the construction
of the Hilbert space presented in sect. 2. They are defined as follows
[\kac]:
 $$
\Theta_{s,p}(a,\tau)=\sum\limits_{\nu\in \longl}\exp\{{{2\pi i\tau s}
 \over {\psi^2}} {(\nu+{p\over
s})}^2+2\pi is(\nu+{p \over s})\cdot a \},
 \eqn\nigran
$$
where $\longl$ stands for
the  long root lattice. These functions are well defined for
 ${\rm Im}\tau >0$, which makes the sum  convergent.
We will consider the case where $p$
belongs to the weight lattice $\wl$.

The Theta functions in \nigran\ satisfy some important properties
\REF\mun{D. Mumford, {\it Tata Lectures on Theta},
Birkh\"auser, Basel, 1983.} [\mun].
The first one, which follows trivially from its definition \nigran,
is the following: a displacement
of $p$ by a vector in $s\longl$ does not change \nigran,
$$
\Theta_{s,p+s\alpha}(a,\tau)=\Theta_{s,p}(a,\tau),\,\,\,\,\,
\alpha\in \longl.
\eqn\axviii
$$
This shows  that  $p$ in
$\Theta_{s,p}(a,\tau)$ lives in
the domain $p\in {\wl \over {s {\longl}}}$.
Another important property is the following. Consider
  $m$ and $n$ two vectors in
 $\longl$, $m,n\in \longl$. Then,
$$
\Theta_{s,p}(a+{2(m+n\tau)\over \psi^2})=
\ex^{2\pi i  s\tau { n\cdot n\over \psi^2}  - 2\pi i  s  n\cdot
a}\Theta_{s,p}(a,\tau).
\eqn\axviiiu
$$
Of particular interest in our analysis are the Weyl antisymmetric combinations
 of Theta functions of level $s$. Let's define them as:
 $$
\Theta^{A}_{s,p}(a,\tau)=\sum\limits_{w\in W} \epsilon
(w)\Theta_{s,w(p)}(a,\tau), \eqn\axix $$
where $\epsilon (w)$ is the signature of the permutation
corresponding the Weyl group element
 $w$. These functions satisfy:
 $$ \Theta_{s,p}^A (a,\tau) =\epsilon
(w)\Theta_{s,w(p)}^A (a,\tau), \eqn\axx
$$
so they are Weyl antisymmetric.
This  property implies some relations
between the antisymmetrized theta functions
of  level $s$. Finally, we recall the behavior of the theta functions under
 modular transformations. The modular group is generated by the
transformation $S$,
$$
\eqalign{ &\tau\rightarrow -{1\over \tau},\cr &a\rightarrow
{a\over\tau},\cr}
\eqn\axxii
$$
and the transformation $T$,
$$
\eqalign{
 &\tau\rightarrow \tau+1,\cr
& a\rightarrow a.\cr}
\eqn\axxiii
$$
The theta functions of level $s$ transform under them as:
 $$
\Theta_{s,p}({a\over\tau},{-1\over\tau})=
\bigg( { { \rm Vol \,\longl^*} \over { \rm Vol\, \longl}} \bigg)^{1 \over
2} \Big({\tau\over is}\Big)^{l\over 2}
\ex^{{i\pi s \over\tau}a^2\psi^2 }\sum\limits_{
q\in {\wl\over {s{\longl}}}}\ex^{-4\pi i{p\cdot q \over
{s \psi^2 }}}
\Theta_{s,q}(a,\tau),
\eqn\axxiv
$$
and,
$$
\Theta_{s,p}(a, \tau+1)=\ex^{2\pi i{p^2\over
s\psi^2}}\Theta_{s,p}(a,\tau).
\eqn\axxv
$$
In \axxiv\  Vol$\, \longl$ is
the volume of the fundamental cell of the long root
lattice $\longl$, and Vol$\, \longl^*$ that of its dual lattice,
$\longl^*$.
The values of their quotient are:
$$
\bigg( { { \rm Vol \,\longl^*} \over { \rm Vol \, \longl}} \bigg)^{1 \over
2} = \cases { N^{-{1 \over 2}}, &$SU(N)$ \cr {1 \over 2}, &$SO(N)$ \cr
2^{-{N \over 4}}, &$Sp(N)$ \cr}\eqn\volumes$$
\endpage

\refout

\end